\documentstyle[12pt
]{article}

\input amssym.def
\input epsf
\newcommand {\SC} {{\Bbb C}}

\newcommand {\SR} {{\Bbb R}}

\newcommand {\SRn} {{\SR^n}}

\newcommand {\SZ} {{\Bbb Z}}

\newcommand {\inR} {\in\SR}

\newcommand {\inRn} {\in\SR^n}

\newcommand {\inZ} {\in\SZ}

\renewcommand {\phi} {{\varphi}}

\newcommand {\e} {{\varepsilon}}
\newcommand {\dt} {{\delta}}
\newcommand {\al} {{\alpha}}
\newcommand {\la} {{\lambda}}

\newcommand {\Scirc} {\raise.2ex\hbox{$\scriptstyle\circ$}}

\newcommand {\ba} {{\boldalpha}}

\newcommand {\bs} {{\bf s}}
\newcommand {\bt} {{\bf t}}

\newcommand {\cD} {{\cal D}}

\newcommand {\cB} {{{\cal B}}}

\newcommand {\cF} {{\cal F}}

\newcommand {\cH} {{\cal H}}

\newcommand {\cL} {{\cal L}}

\newcommand {\cS} {{\cal S}}

\newcommand {\tc} {{\tilde c}}

\newcommand {\tB} {{\widetilde B}}

\newcommand {\tJ} {{\widetilde J}}

\newcommand {\txi} {{\tilde\xi}}
\newcommand {\tpsi} {{\tilde\psi}}

\newcommand {\om} {{\omega}}

\newcommand {\hf} {{\hat f}}

\newcommand {\Ga} {{ G^{\al}}}

\newcommand {\hg} {{\hat g}}

\newcommand {\hchi} {{\widehat\chi}}

\newcommand {\hpsi} {{\widehat\psi}}
\newcommand {\htpsi} {\,\widehat{\!\tpsi}}

\newcommand {\hphi} {{\widehat\phi}}

\newcommand {\re} {{\Re\!\mbox{\small\it e}\,}}

\newcommand {\Dquad} {\qquad\qquad}
\newcommand {\Tquad} {\Dquad\qquad}
\newcommand {\Qquad} {\Dquad\Dquad}
\newcommand {\Pquad} {\Qquad\qquad}

\newcommand {\mand} {{\quad\mbox{and}\quad}}

\renewcommand {\liminf}{\mathop{\underline{\mbox{\rm lim}}}}

\newcommand {\tfrac}[2]{{\textstyle\frac{#1}{#2}}}

\newcommand {\supp} {\mbox{Supp }}
\newcommand {\rank} {\mbox{rank }}
\newcommand {\tr} {\mbox{tr}}
\newcommand {\meas} {\mbox{meas }}

\newcommand {\diag} {\mbox{diag }}

\newcommand {\Vol} {\mbox{Vol }}

\newcommand {\Proof} {\noindent{\bf P{\footnotesize\bf ROOF}: } \ }
\newcommand {\Proofof}[1] {\noindent{\bf P{\footnotesize\bf ROOF} of {#1}: } \ }
\newcommand {\ProofEnd} {
             \begin{flushright} \vskip -0.2in $\Box$ \end{flushright}}

\renewcommand {\mid} {{\,\,\,\colon\,\,\,}}

\newcommand {\Ol} {\overline}

\newcommand {\Ds} {\displaystyle}
\newcommand {\Ts} {\textstyle}

\newcommand{\Ba}[1]{\begin{array}{#1}}
\newcommand{\Ea}{\end{array}}
\newcommand{\Be}{\begin{equation}}
\newcommand{\Ee}{\end{equation}}
\newcommand{\Bea}{\begin{eqnarray}}
\newcommand{\Eea}{\end{eqnarray}}
\newcommand{\Beas}{\begin{eqnarray*}}
\newcommand{\Eeas}{\end{eqnarray*}}
\newcommand{\Benu}{\begin{enumerate}}
\newcommand{\Eenu}{\end{enumerate}}
\newcommand{\Bi}{\begin{itemize}}
\newcommand{\Ei}{\end{itemize}}

\newcommand{\Bv}{\begin{verse}}
\newcommand{\Ev}{\end{verse}}

\newcommand{\BR}{\begin{Remark} \em}
\newcommand{\ER}{\end{Remark}}

\newcommand{\BE}{\begin{example} \em}
\newcommand{\EE}{\end{example}}

\newcounter{remark}

\newtheorem{theorem}[equation]{T{\hskip 0pt\footnotesize\bf HEOREM}}
\newtheorem{proposition}[equation]{P{\hskip 0pt\footnotesize\bf ROPOSITION}}
\newtheorem{corollary}[equation]{C{\hskip 0pt\footnotesize\bf OROLLARY}}
\newtheorem{lemma}[equation]{L{\hskip 0pt\footnotesize\bf EMMA}}
\newtheorem{Remark}[equation]{R{\hskip 0pt\footnotesize\bf EMARK}}
\newtheorem{definition}[equation]{D{\hskip 0pt\footnotesize\bf EFINITION}}
\newtheorem{example}[equation]{E{\hskip 0pt\footnotesize\bf XAMPLE}}

\newcommand{\La}{{\Lambda}}
\newcommand{\xiljk}{{\xi^\ell_{j,k}}}
\newcommand{\hh}{{\hat{h}}}

\newcommand{\sh}{{\mbox{sh }}}

\newcommand{\fnr}{\frac{n}r}

\newcommand{\Ad}{A_\nu^2}
\newcommand{\Ld}{L_\nu^2}

\newcommand{\Apq}{{A^{p,q}_\nu}}

\newcommand{\Bpq}{{B^{p,q}_\nu}}
\newcommand{\cBpq}{{{\cB}^{p,q}_\nu}}

\newcommand{\Bp}{{B^{p',q'}_{-\nu{q'}/q}}}

\newcommand{\Lpq}{{L^{p,q}_\nu}}

\newcommand{\cLt}{\widetilde{\cal L}}
\newcommand{\cLL}{{\cal L}}
\renewcommand{\cL}{{\cal E}}
\renewcommand{\cLt}{\widetilde{\cal E}}

\newtheorem{notation}[equation]{N{\hskip 0pt\footnotesize\bf OTATION}}

\def\lan#1#2{{\langle\,{#1}\,,\,{#2}\,\rangle}}

\newcommand{\bO}{\partial\Omega}
\newcommand{\clO}{{\Ol{\Omega}}}

\renewcommand{\O}{\Omega}
\renewcommand{\Ga}{\Gamma}
\newcommand{\Tu}{T_\O}
\newcommand{\So}{{\cS_\O}}
\newcommand{\Do}{{\cD_\O}}
\newcommand{\Sop}{{\cS'_{\bar{\O}}}}
\newcommand{\Sbop}{{\cS'_{\bO}}}

\newcommand{\be}{{\bf e}}
\renewcommand{\ba}{{\bf \al}}
\newcommand{{\tBox}}{{\widetilde{\raisebox{-0.2ex}[1.25ex][0ex]{$\Box$}}}}

\newcommand{\Dt}{\Delta}
\newcommand{\dmes}{\Dt(y)^{\nu-\fnr}\,dy}
\newcommand{\dxi}{{\frac{d\xi}{\Dt(\xi)^{\frac{n}r}}}}
\newcommand{\ind}{{\Ts{\frac{n}r-1}}}
\newcommand{\tqn}{\tilde{q}_{\nu,p}}

\newcommand{\lesssim}{\;\Subb{<}{\sim}\;}
\newcommand{\gtrsim}{\;\Subb{>}{\sim}\;}
\newcommand {\Subb}[2] {\raisebox{-.8 ex}
{$\stackrel{\displaystyle #1} {\scriptstyle #2}$}}

\begin{document}

\title{Littlewood-Paley decompositions and Besov
spaces related to symmetric  cones}
\author{D. B\'ekoll\'e\quad A. Bonami$^*$\quad
G. Garrig\'os\footnote{ Research supported by the European
Commission, within the Networks ``TMR: Harmonic Analysis
1998-2002'' and ``IHP: HARP 2002-2006''. Third author supported by
\emph{Programa Ram\'on y Cajal} and grant ``BMF2001-0189'', MCyT
(Spain).}
 \quad F. Ricci$^*$}

\maketitle


\begin{abstract}
Starting from a Whitney decomposition of a symmetric cone
$\Omega$, analog to the dyadic partition $[2^j, 2^{j+1})$
 of the positive real line, in this paper we  develop
an adapted Littlewood-Paley theory for functions with spectrum in
$\Omega$. In particular, we define a natural class of Besov spaces
of such functions, $\Bpq$, where the role of usual derivation is
now played by the generalized wave operator of the cone
$\Dt(\frac{\partial}{\partial x})$. Our main result shows that
$\Bpq$ consists precisely of the distributional boundary values of
holomorphic functions in the Bergman space $\Apq(\Tu)$, at least
in a ``good range'' of indices
 $1\leq q<q_{\nu,p}$. We obtain the sharp $q_{\nu,p}$ when $p\leq 2$,
and conjecture a critical index for $p>2$. Moreover, we show the
equivalence of this problem with the boundedness of Bergman
projectors $P_\nu\colon\Lpq\to\Apq$, for which our result implies
a positive answer when $q_{\nu,p}'<q<q_{\nu,p}$. This extends to
general cones previous work of the authors in the light-cone.
Finally, we conclude the paper with a finer analysis in
light-cones, for which we establish a link between our conjecture
and the cone multiplier problem. Moreover, using recent work by
Tao, Vargas and Wolff, we improve in dimension 3 the range of
$q$'s for which the Bergman projection is bounded.
\end{abstract}

\section{Introduction}
\setcounter{equation}{0} \setcounter{footnote}{0}
\setcounter{figure}{0}

Let $\O$ be an irreducible symmetric cone in a Euclidean vector
space $V$ of dimension $n$, endowed with an inner product
$(\,\cdot\,|\,\cdot\,)$ for which the cone $\O$ is self-dual. We
can identify $V$ with $\SR^n$, by endowing the latter with such
inner product. We denote by $\Tu=V+i\O$ the corresponding tube
domain in the complexification of $V$, which we may also identify
with $\SC^n$. As in the text \cite{FK}, we shall write the rank
and determinant associated with a cone by
\[
 r=\rank{\O},\mand\Dt(x)=\det{x}, \quad x\in V.
\]

Two examples of the above situation are the light-cones and the
cones of positive definite symmetric matrices. The first ones are
defined in $\SRn$, for $n\geq3$, by
\[
\La_n=\{y=(y_1,y')\in\SRn\mid y_1^2-|y'|^2>0,\;y_1>0\}.
\]
These are symmetric cones of rank $2$ with determinant given by
the \emph{Lorentz form} $\Dt(y)=y_1^2-|y'|^2$. The cones
$\mbox{Sym}_+\,(r,\SR)$ of positive definite symmetric matrices,
have rank $r$ and the usual determinant for matrices. In this last
case, the underlying vector space $V$ is the space of symmetric
matrices $\mbox{Sym}\,(r,\SR)$, with dimension $n=\frac{r(r+1)}2$,
and with a Euclidean norm defined by the Hilbert-Schmidt inner
product. We observe that this does not coincide with the canonical
inner product in the usual identification between $V$ and $\SR^n$.
\medskip

 The goal
of this paper is to present, in the general setting of symmetric
cones, a special Littlewood-Paley decomposition adapted to the
geometry of $\O$. From this we shall obtain new results in
analytic problems, such as the boundedness of Bergman projectors
and the characterization of boundary values for Bergman spaces in
the tube domain $\Tu$. To describe our setting, we shall denote by
$\So$ the space of Schwartz functions $f\in\cS(\SRn)$ with
$\supp\hf\subset\clO$, and normalize the Fourier transform by
\[ \hf(\xi)=\cF f(\xi)=
\frac1{(2\pi)^n}\,\int_\SRn e^{-i(x|\xi)} f(x)\,dx,\quad\xi\inRn.
\]
Our key tool will be the following special decomposition for
functions in $\So$ \Be f=\sum_j f*\psi_j, \quad\forall\;f\in\So,
 \label{dec}
\Ee where $\hpsi_j$ are test functions supported on ``frequency
blocks''  $B_j$, and these form a suitable \emph{Whitney covering}
on the cone $\O$. In analogy with the dyadic decomposition of the
half-line $(0,\infty)$ (i.e., the $1$-dimensional cone), the sets
$B_j$ are constructed from the homogeneous structure of $\O$, as
the ``balls'' $B_j=\{\xi\in\O\mid d(\xi,\xi_j)<1\}$, obtained from
a $G$-invariant distance $d$ and  a $d$-lattice $\{\xi_j\}$  in
$\O$. These will turn out to be the right sets for the
discretization of many operators in the cone, since functions
which appear in their multiplier expressions, such as
 $\Dt(\xi)$ or $(\xi|y)$ (for fixed $y\in \O$), remain essentially constant
 when $\xi\in B_j$.

A characteristic example of this situation is the  {\it
generalized wave operator} on the cone: $\Box=\Dt(\frac 1i
\frac{\partial}{\partial x})$, which is the differential operator
of degree $r$ defined by the equality: \Be \Dt{\Ts\left(\frac 1i
\frac{\partial}{\partial x}\right)}
[e^{i(x|\xi)}]=\Dt(\xi)e^{i(x|\xi)}, \quad\xi\inRn. \label{box}
\Ee This corresponds, in cones of rank 1 and 2, to
\[
\Box=\frac1i\,\frac{d}{dx}\quad \mbox{in}\;(0,\infty),
\quad\mbox{and}\quad
\Box=-\frac14\,\left(\frac{\partial^2}{\partial x^2_1}-
\frac{\partial^2}{\partial x^2_2}-\ldots-
\frac{\partial^2}{\partial x^2_n}\right)\quad \mbox{in}\;\La_n.
\]
The Littlewood-Paley decomposition (\ref{dec})
 provides a formal ``discretization'' of the action of $\Box$
on functions with spectrum in $\O$:
\[
\Box f=\cF^{-1}(\Dt(\xi)\hf(\xi))= \sum_j\Dt(\xi_j)\,
f*\psi_j*\cF^{-1}(m_j),\quad f\in\So,
\]
where $\{m_j\}$ is a uniformly bounded family of multipliers.

\medskip

From these facts it is natural to introduce a new family of
Besov-type spaces, $\Bpq$, adapted to the Littlewood-Paley
decomposition (\ref{dec}). These are defined as the equivalence
classes of tempered distributions, which have finite seminorms \Be
 \|f\|_{\Bpq}=
 \left[\sum_j\Dt^{-\nu}(\xi_j)\,\|f*\psi_j\|_p^q
\right]^\frac1q.
 \label{bpqnorm}
 \Ee
Our first result shows that these spaces satisfy analogous
properties to the one dimensional homogeneous Besov spaces, with
the role of usual derivation played by the wave operator $\Box$.
We warn the reader that, for convenience in the applications that
follow, we are using a non-standard normalization of indices in
our definition of $\|\cdot\|_\Bpq$ (compared, e.g., with
\cite{FJW}).

\begin{theorem}
\label{thb} Let $\nu\inR$ and $1\leq p,q<\infty$. Then \Benu \item
$\Bpq$ is a Banach space and does not depend on the choice of
$\{\psi_j\}$ and $\{\xi_j\}$. \item $\Box:\Bpq\rightarrow
B^{p,q}_{\nu+q}$ is an isomorphism of Banach spaces.

\item If $p,q>1$, then $(B^{p,q}_{\nu})^*$ is isomorphic to
 $B^{p',q'}_{-\nu q'/q}$ with the usual duality pairing.
\Eenu
\end{theorem}


The rest of the paper is devoted to applications of this theory to
two open problems involving the class of \emph{Bergman spaces}. In
this paper, a weighted mixed-norm version of these spaces is
defined by  the integrability condition:
 \Be
 \|F\|_\Lpq:=
 \left[\int_\O\Bigl(\int_\SRn|F(x+iy)|^p\,dx
 \Bigr)^\frac{q}{p}\,
 \Dt^{\nu-\fnr}(y)\,dy\right]^\frac1q<\infty.
 \label{lpqnorm}
 \Ee
Thus, when $1\leq p,q<\infty$ and $\nu\in\SR$, we denote by
$\Apq(\Tu)$ the closed subspace of $\Lpq$ consisting of
holomorphic functions in the tube $\Tu$. We observe that these
spaces are non null only when $\nu>\ind$ (see, e.g.,
\cite{BBGNPR}). The usual $A^p$ space corresponds to $p=q$ and
$\nu=\frac nr$. To simplify notation we shall write
$A^p_\nu=A^{p,p}_\nu$, and similarly  $L^p_\nu= L^{p,p}_\nu$.

Two main questions concerning these spaces will be studied here:
\Benu \item The characterization of boundary values of functions
in $\Apq$, as distributions in the Besov spaces $\Bpq$.

\item The existence of bounded extensions into $\Lpq$ spaces for
the (weighted) Bergman projector, that is, the orthogonal
projector $P_\nu\colon L^2_\nu\to A^2_\nu$. \Eenu

Regarding the first question, it has been known for some decades
the relation, in the 1-dimensional setting, between boundary
values of Bergman functions and homogeneous Besov spaces (see e.g.
\cite{RT,BUI}, or the lecture notes \cite{BBGNPR}). To see this in
higher dimensions, and restricted to tube domains over cones
$\Tu$, one writes a homolomorphic function $F\in \Apq$ in terms of
its {\it Fourier-Laplace transform}: \Be F(z)=\cLL g(z)=\int_\O
e^{i\,(z|\xi)}g(\xi)\,d\xi, \quad z\in\Tu, \label{Ft} \Ee for some
distribution $g$ supported in $\clO$. Observe that the new
distribution $f=\cF^{-1}g$ plays the role of a ``Shilov boundary
value'' for $F$, and hence it is a natural candidate to belong to
$\Bpq$. Now, by Theorem \ref{thb} this is equivalent to
$\cF^{-1}(e^{-(y|\cdot)}\chi_\O)$ having a finite $B^{p',q'}_{-\nu
q'/q}$ norm, which by explicit computation can only happen when
$q$ is below a certain critical index
\[ \tqn=\frac{\nu+\fnr-1}{(\fnr\frac1{p'}-1)_+}
\]
(with $\tqn=\infty$, if $\fnr\leq p'$). This constitutes our first
result, whose detailed justification will be presented in sections
3.4 and 4.1.

%
%

 \begin{theorem}
 \label{th1}
 Let $\nu>\ind$, $1\leq p<\infty$ and $1\leq q<\tqn$.
 Then, for every $F\in\Apq$ there exists
a (unique) tempered distribution $f\in\cS'(\SRn)$ such that
$f=\sum_jf*\psi_j$ in $\cS'(\SR^n)$, $\|f\|_\Bpq<\infty$ and
$F=\cLL \hf$. Moreover
 we have
 \Beas
1. & & \lim_{{y\to0}\atop{y\in\O}}F(\cdot+iy)= f,\quad\mbox{both
in}\; \cS'(\SR^n)\mbox{ and }\; \Bpq;\\
2. & & \|f\|_{\Bpq}\leq C\,\|F\|_{\Apq},\quad\mbox{for all}\quad
F\in\Apq. \Eeas
 \end{theorem}

The converse result is more interesting, and turns out to be
equivalent to the second of the questions posed above. We only
have a partial answer, for which we need to introduce two new
critical indices
\[
q_\nu=\frac{\nu+\fnr-1}{\fnr-1},\quad
q_{\nu,p}=\min\{p,p'\}\,q_\nu\,.\] Observe that in 1-dimension the
three indices are equal to $\infty$, while in general we have the
ordering
\[
2<q_\nu\leq q_{\nu,p}\leq \tqn\,.
\]
The role of these new indices will be clarified later in relation
with the Bergman projectors.

\begin{theorem}
\label{th2} Let $\nu>\ind$, $1\leq p<\infty$ and $1\leq
q<q_{\nu,p}$. Given a distribution $f\in\cS'(\SR^n)$ such that
$f=\sum_jf*\psi_j$ and $\|f\|_\Bpq<\infty$, then the holomorphic
function $F=\cLL \hf$ belongs to $\Apq$, and moreover, there
exists a constant $C>0$ so that
\[
\frac 1C\,\|f\|_\Bpq\,\leq \,\|\cLL\hf\|_\Apq\,\leq\,C\,
 \|f\|_\Bpq,\quad f\in\Bpq.
\]
\end{theorem}

This theorem is sharp for $1\leq p\leq2$, in the sense that for
each $q\geq q_{\nu,p}=pq_\nu$ there is a distribution with
$\|f\|_\Bpq<\infty$ and $\|\cLL \hf\|_\Lpq=\infty$. We shall
present these examples in $\S4.4$. When $p>2$ we will construct
similar examples, but only for values of $q\geq q_{\nu,2}=2q_\nu$,
leaving open the question when $p'q_\nu\leq q<\min\{2q_\nu,
\tqn\}$. New positive results in dimension 3 can be obtained using
restriction theorems, and will be presented in the Appendix.

\medskip

Finally, we turn to the second application of our theory, the
boundedness of Bergman projectors in $\Lpq$. This is a challenging
question which has been open for many years, and which still is
not completely solved. The three indices defined above correspond
to three steps of difficulty for this question. For instance, a
trivial counterexample shows that $P_\nu$ can only be bounded in
$\Lpq$ for $\tqn'<q<\tqn$. This follows from
the fact that
the Bergman kernel belongs to $L^{p',q'}_\nu(\Tu)$ only when
$q<\tqn$ (see $\S4.3$ below). From the other two indices, the
smallest one gives the natural range $q'_\nu<q<q_\nu$ for
boundedness of the \emph{positive operator} $P^+_\nu$, obtained by
replacing the Bergman kernel with its absolute value
$|B_\nu(z,w)|$ \cite{BB1,BT}. Finally, $q_{\nu,p}$ appears as the
interpolation index between $q_\nu$ and $q_{\nu,2}=2q_\nu$, giving
the latter the sharp range of boundedness in the spaces
$L^{2,q}_\nu$ \cite{BBPR}. The results in this last paper, which
combine the Plancherel identity with a suitable ``discretization
of multipliers'' in light-cones, have been the germ of the
Littlewood-Paley decomposition we are introducing here. Our main
contribution to this problem is, besides an extension to general
symmetric cones, a direct formulation in terms  of
Littlewood-Paley inequalities, which allows further improvements
as those considered in the Appendix. We gather these results in
our next theorem, which may be stated only for $q\geq2$ by
self-adjointness of $P_\nu$.

\begin{theorem}
\label{th3} Let $\nu>\ind$, $1\leq p<\infty$ and $2\leq
q<q_{\nu,p}$. Then, the inequality \[ \|\cLL\hf\|_\Lpq\leq
C\,\|f\|_\Bpq,\quad f\in\So,
\] holds true if and only if $P_\nu$ can be
boundedly extended from $\Lpq$ onto $\Apq$. In particular, $P_\nu$
is bounded in $\Lpq$ for all $1\leq p<\infty$ and
$q_{\nu,p}'<q<q_{\nu,p}$. Moreover, $P_\nu$ does not admit bounded
extensions to $\Lpq$ when: \Benu \item $1\leq p\leq 2$ and $q\geq
q_{\nu,p}$;

\item $2<p<\infty$ and $q\geq \min\{2 q_\nu,\tqn\}.$ \Eenu
\end{theorem}

\begin{figure}[ht]
   \centerline{
\epsfysize=2.5in \epsfxsize=2.5in \epsfbox{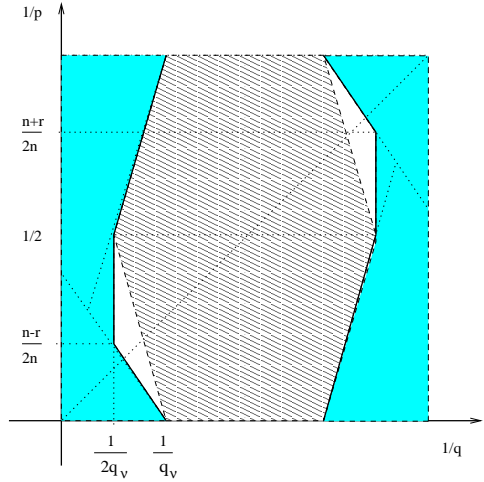}}
\caption{Region of boundedness of $P_\nu$.} \label{fig}
\end{figure}

Figure \ref{fig} illustrates the regions of boundedness,
unboundedness, and the open gap where for the moment no answer is
known (compare with \cite{BBPR}). In the particular case of
$L^p_\nu$-spaces ($p=q$) the gap becomes $1+ q_\nu\leq p<\min\{2
q_\nu, q_\nu+\frac{n}{n-r}\}$, together with the corresponding
dual interval.
We point out that our results produce new counterexamples compared
with \cite{BBPR}, first by considering the limit cases when $1\leq
p\leq 2$ and $q\geq2$, and second by removing the region with
$q\geq q_{\nu,2}$ when $p>2$. The counterexamples presented here
are besides valid for general symmetric cones.

\medskip

We do not wish to conclude this introductory section without
mentioning our approach to the open question when $p,q>2$. Our proofs
provide sufficient conditions which are variants of the elementary inequality
\Be
 (\sum_j\|f*\psi_j\|_p^s)^\frac1s\leq\,C\,\|f\|_p,
 \quad f\in\cS(\SR^n)
 \label{LP}
 \Ee
for indices $2\leq s\leq p$. The inequality holds trivially for $s=p$ (by
interpolation between $s=\infty$ and $s=p=2$), while any improvement
in $s$ smaller than $p$ will directly imply
boundedness of the Bergman projection outside the hexagonal region.
In particular, going down
to $s=2$ will fill the gap in Figure \ref{fig} up the vertical
line $q=2q_\nu$. Variants of such inequalities (typically with
$s=2$) have been widely studied in light-cones, as well as their
analogs in spheres. In particular, square-function versions of
it:\Be
 \|\sum_j f*\psi_j\|_{p'}\leq \,C\,
\Bigl\|\,\bigl(\sum_j |f*\psi_j|^2\bigr)^\frac12\Bigr\|_{p'},\quad
f\in\So,
 \label{LPs}\Ee
are intimately related with the cone multiplier problem and the
Bochner-Riesz means. We shall show in the Appendix how such
results lead to a slight improvement of our theorems in the
3-dimensional light-cone. Morally, this reduces the problem of
 boundedness of Bergman projectors to  very challenging
 questions related to inequalities like (\ref {LP}) or (\ref{LPs}), where the
 complex analysis formalism has been completely removed.
\smallskip

Finally, we conclude by mentioning that a simplified version of our results,
specialized to the case when $p=2$ (for which the use of
Plancherel is available), has been published separately in
\cite{BBG}. We also refer to the survey paper \cite{B} for
complementary information concerning the critical indices and the
relation with Hardy-type inequalities.

\section{Whitney decompositions on the cone}
\setcounter{equation}{0} \setcounter{footnote}{0}
\setcounter{figure}{0}

In this section we introduce the notation and a list of technical
results on symmetric cones, mostly taken from the text \cite{FK}.
We also give a detailed construction of the ``Whitney
decomposition'' adapted to the analysis of the problems stated
above. The main lines and applications of such constructions
appear in previous papers: see \cite{BBPR} for the light-cone, and
\cite{BBG} for general symmetric cones.

\subsection{Background on symmetric cones}

Let $\O$ be a fixed symmetric cone in a real Euclidean vector
space $V$, endowed with the inner product $(\,\cdot\,|\,\cdot\,)$.
That is, $\O$  is a homogeneous open convex cone which is
self-adjoint with respect to $(\,\cdot\,|\,\cdot\,)$. Let $G(\O)$
be the group of linear transformations of the cone, and $G$ its
identity component. By definition, $G(\O)$ acts transitively on
$\O$. Further, it is well known that
 there is a solvable subgroup $T$ of $G$
acting \emph{simply transitively on $\O$}. That is, every $y\in\O$
can be written uniquely as $y=t\be$, with $t\in T$ and a fixed
$\be\in\O$. This gives an identification $\O\equiv T=G/K$, where
$K$ is a maximal compact subgroup of $G$. Moreover, $K=\{g\in
G\mid g\be=\be\}=G\cap O(V)$.
 All these properties can be found in the first chapter of the text \cite{FK}.

It is well-known that for every symmetric cone $\O$, its
underlying vector space $V$ can be endowed with a multiplication
rule which makes it a \emph{Euclidean Jordan algebra} with
identity element $\be$. With such multiplication, $\clO$ coincides
with the set $\{x^2\mid x\in V\}$ of all squares in $V$. We may
assume that the inner product in $V$ is given by
$(x|y)=(xy|\be)=\tr(xy)$ (see \cite[Ch. III]{FK}). The reader less
familiar with these concepts can think on the example of
positive-definite symmetric matrices, which we present in
\ref{p-d} below.

Suppose now that the cone is irreducible, has rank $r$ and its
underlying space has dimension $n$. A precise form for the group
$T$ can be obtained from the Jordan algebra structure of $V$.
Following \cite[Ch. VI]{FK}, we let $\{c_1,\ldots,c_r\}$ denote a
fixed Jordan frame in $V$, and $V=\oplus_{1\leq i\leq j\leq
r}V_{i,j}$ its associated Peirce decomposition. Then $T$ may be
taken as the corresponding solvable Lie group, which factors as
the semidirect product $T=NA=AN$ of a nilpotent subgroup $N$ (of
lower triangular matrices), and an abelian subgroup $A$ (of
diagonal matrices). The latter takes the explicit form
\[
A=\{P(a)\mid a=\sum_{i=1}^ra_ic_i,\,\,\,a_i>0\},
\]
where  $P$ is the quadratic representation of $V$. This leads as
well to the classical decompositions of the semisimple Lie group
$G=NAK$ and $G=KAK$.

Still following \cite[Ch. VI]{FK}, we shall denote  by
$\Dt_1(x),\ldots,\Dt_r(x)$ the
 {\it principal minors} of $x\in V$, with
respect to the fixed Jordan frame $\{c_1,\ldots,c_r\}$.
These are invariant functions under the group $N$:
\[
\Dt_k(nx)=\Dt_k(x), \quad n\in N, \quad x\in V,\quad k=1,\ldots,r,
\]
and satisfy a homogeneity relation  under $A$
\[
\Dt_k(P(a)x)=a_1^2\cdots a_k^2\,\Dt_k(x),\quad \mbox{if}\quad
a=a_1c_1+\ldots+a_rc_r.
\]
The determinant function $\Dt(y)=\Dt_r(y)$ is also invariant under
$K$, and moreover, satisfies the formula \Be\label{inv-det}
\Dt(gy)=\Dt(g\be)\Dt(y)=\mbox{Det(g)}^{\frac rn} \Dt(y).\Ee It
follows from this formula that an invariant measure in $\O$ is
given by $\Dt(y)^{-\frac{n}r}\,dy$. Finally, we recall a version
of Sylvester's Theorem for symmetric cones, which allows to write
these as:
\[
\O=\{x\in V\mid \Dt_k(x)>0,\,\,\,k=1,\ldots,r\}.
\]

\medskip

\begin{example}\label{p-d}{\bf The cone of positive-definite symmetric
matrices} \smallskip

 {\rm We describe the above concepts for the
cone $\O= \mbox{Sym}_+ (r, {\bf R})$, contained in the vector
space $V= \mbox{Sym} (r,
 {\bf R})$.  The \emph{Jordan algebra structure} in
$V$ corresponds to the symmetric product $X\circ
Y=\frac12(XY+YX)$, with the usual identity matrix $\be=I$. A
standard \emph{Jordan frame} is the set $D_j$ of diagonal matrices
all of whose entries are $0$ except for the $j$-th one equal to
$1$. The Peirce decomposition in $V$ is just the decomposition of
a symmetric matrix in terms its $(i,j)$ entries.

In this example, the \emph{automorphism group} $G(\O)$ can be
identified with $\mbox{Gl}(r,\SR)$ via the adjoint action: \Be
g\in{\mbox{Gl}}(r,\SR),Y\in\mbox{Sym}(r,\SR)\longmapsto g\cdot
Y=gYg^*\in\mbox{Sym}(r,\SR). \Ee Then, the group $T$ consists in
the lower triangular matrices in $\mbox{Gl}(r,\SR)$, and the
factorization $Y=t\cdot I$ is precisely the Gauss decomposition of
a symmetric matrix.  The subgroup $N$ consists of all triangular
matrices in $\mbox{Gl}(r,\SR)$ with $1$'s on the diagonal, while
$A$ is given by the diagonal matrices
$P(a)=\diag\{a_1,\ldots,a_r\}$. Finally, the associated
\emph{principal minors} are the usual principal minors from linear
algebra, that is, the determinants of the $k\times k$ symmetric
matrices obtained by restriction to the first $k$ coordinates. One
verifies easily with this example the homogeneity properties with
respect to $N$ and $A$ stated above.}
\end{example}

\subsection{The invariant metric and the covering lemma}

With the identification  $\O\equiv G/K$, the cone can be regarded
as a Riemannian manifold with the $G$-invariant metric defined by
\[
\langle\xi,\eta\rangle_y:=(t^{-1}\xi|t^{-1}\eta)
\]
if $y=t\be$ and $\xi,\eta$ are tangent vectors at $y\in\O$. We
shall denote by $d(\,\cdot\,,\,\cdot\,)$ the corresponding
distance, and by $B_\dt(\xi)$ the ball centered at $\xi$ of radius
$\dt$. Note that, for each $g\in G$, the invariance implies
$B_\dt(g\xi)=gB_\dt(\xi)$.

We shall need some weak local invariance properties of the
quantities that we have defined on the cone. One consequence is
the possibility of obtaining a Whitney-type decomposition for
general symmetric cones in terms of invariant balls. Part of this
material  was already presented in \cite{BBG}.

\begin{lemma}
\label{l3} Let $\dt>0$. Then there is a constant
$\gamma=\gamma(\dt,\O)>0$ such that
\[
\mbox{if}\quad d(\xi,\xi')\leq\dt\quad \Rightarrow\quad
\frac1\gamma\leq\frac{\Dt_k(\xi)}{\Dt_k(\xi')}\leq\gamma,\quad
k=1,\ldots,r.
\]
\end{lemma}
\Proof By invariance of the metric and the forms $\Dt_k$ under
$N$, we may assume $\xi'=P(a)\be$. Further, since
\[
\frac{\Dt_k(\xi)}{\Dt_k(P(a)\be)}=\frac{\Dt_k(P(a)^{-1}\xi)}{\Dt_k(\be)},
\]
we may even assume $\xi'=\be$. Now, the estimations above and
below for $\Dt_k$ in a ball ${\Ol{B}}_{\dt}(\be)$ follow easily
from the continuity of $\xi\mapsto\Dt_k(\xi)$, and a compactness
argument. \ProofEnd

The next lemma states the local equivalence between two Riemannian
metrics. The proof follows from standard arguments (see, e.g.,
\cite[9--22]{spivak}).
\begin{lemma}
\label{l7} Let $\dt_0>0$ be fixed. Then, there exist two constants
$\eta_1>\eta_0>0$,
 depending only on $\dt_0$ and $\O$,
 so that for every $0<\dt\leq\dt_0$ we have
\[
\bigl\{|\xi-\be|<\eta_0\dt\bigr\} \,\subset\,
B_\dt(\be)\,\subset\, \bigl\{|\xi-\be|<\eta_1\dt\bigr\}.
\]
\end{lemma}

We can now estimate the volume of an invariant ball. Recall that
the invariant measure in $\O$ is given by
\[
\meas(B)=\int_B\Dt(\xi)^{-\frac{n}r}\,d\xi, \quad
B\subset\O\,\,\,\mbox{measurable}.
\]
Therefore, from the previous results it follows that,  for all
$y\in\O$ and $0<\dt\leq\dt_0$,
\[
\meas(B_\dt(y))=\meas(B_\dt(\be))\sim\Vol(B_\dt(\be))\sim\dt^n,
\]
where the equivalences denoted by ``$\sim$'' are modulo constants
depending only on $\O$ and the fixed number $\dt_0$. Observe,
however, that this estimate cannot hold uniformly in
$\dt_0>\!\!>1$, since the invariant measure is in general not
doubling. We can now prove a covering lemma which will be of
crucial importance for the rest of the paper.

\begin{lemma}
\label{l1}\,\,
 {\bf:\,\, Whitney Decomposition of the cone.}\,
Let $\dt>0$ and $R\geq2$. Then, there exist sequences of points
$\{\xi_j\}_j$ in $\O$ such that \Benu \item[(i)\,]
$\{B_\dt(\xi_j)\}_j$ is a disjoint family in $\O$;

\item[(ii)] $\{B_{R\dt}(\xi_j)\}_j$ is a covering of the cone
$\O$. \Eenu
 Moreover, for each such sequence the balls
$\{B_{R\dt}(\xi_j)\}_j$ have the finite intersection property.
That is, if $\dt, R\leq R_0$, then there exists an integer
$N=N(R_0, \O)$ so that at most $N$ of these balls can intersect an
arbitrary set $E\subset\O$ with diameter \Be \mbox{diam
}(E)=\sup\{d(\xi,\eta)\mid\xi,\eta\in E\}\leq R_0\dt.
\label{diam1}\Ee
%
\end{lemma}

\Proof Consider $\{\xi_j\}_j$ a maximal subset of $\O$ (under
inclusion) among those with the property that their elements are
distant at least $2\dt$ from one another. Let us denote $B'_{j}$
the balls $B_{\dt}(\xi_j)$. They are pairwise disjoint, while, by
maximality, the  balls $\{B_j=B_{2\dt}(\xi_j)\}_j$  cover $\O$.
Note also that, necessarily, the set $\{\xi_j\}_j$ is countable.

For the finite overlapping, let $E$ be a set as in (\ref{diam1}).
Denote by $J$ the set of indices $\{j\mid B_{R\dt}(\xi_j)\cap
E\not=\emptyset\}$, and fix a point $\xi\in
B_{R\dt}(\xi_{j_0})\cap E$ for some $j_0\in J$. Then, the
condition on the diameter gives
\[
\bigcup_{j\in J}B_\dt(\xi_j) \,\subset\,B_{(2R_0+1)\dt}(\xi)\,.
\]
Now, by disjointness and invariance of the measure we have
\[
|J|\,\meas(B_{\dt}(\be))=\meas(\cup_{j\in J}B'_j)\leq
\meas(B_{(2R_0+1)\dt}(\xi))=\meas(B_{(2R_0+1)\dt}(\be)).
\]
Thus, the remarks preceding the lemma give us a bound for $N$
depending only on $\O$ and $R_0$.

\ProofEnd

\BR \Benu

\item A sequence of points $\{\xi_j\}_j$  with the above
properties will be called a {\it $(\dt,R)$-lattice of the cone}.
Observe that one can always define an associated partition by
letting
\[
E_1=B_1,\quad\ldots\quad,E_j=B_j\setminus E_{j-1},\quad\ldots
\]
We shall call $\{E_j\}_j$ a {\it Whitney decomposition of $\O$}.

\item If $\{\xi_j\}_j$ is a $(\dt,R)$-lattice, then so is
$\{\xi_j^{-1}\}_j$. Indeed, this follows from the fact that
$y\mapsto y^{-1}$ is an isometry of the cone (see Chapter III of
\cite{FK}). Therefore,
 $B_\dt(\xi_j^{-1})=B_\dt(\xi_j)^{-1}$, and the
conditions of Lemma \ref{l1} hold.

\item One can look at the sequences
 $\{\xi_j\}_j$ and $\{\xi_j^{-1}\}_j$ as
a {\it couple of  dual lattices}. In fact,
$(\xi_{j}\vert\xi_{j}^{-1})=r$, while using
$\Dt(y^{-1})=\Dt(y)^{-1}$ we also have
$\Vol(B_j)\sim\Dt(\xi_j)^{\fnr}$ and
$\Vol(B_j^{-1})\sim\Dt(\xi_j)^{-\fnr}$. Moreover, from the next
lemma it will follow that actually $(\xi\vert y)\sim 1$ when  $\xi
\in B_{j}$ and $y\in B_{j}^{-1}$. \Eenu \ER


\begin{lemma}
\label{l6} Let $\dt>0$. There exists $\gamma=\gamma(\O,\dt)>0$
such that, for $y\in{\Ol{\O}}$ and $\xi,\xi'\in\O$ with
$d(\xi,\xi')\leq\dt$, then \Be
\frac1\gamma\,\leq\,\frac{(\xi|y)}{(\xi'|y)}\, \leq\gamma.
\label{inner} \Ee In particular,
$\frac1\gamma\,\leq\,|\xi|/|\xi'|\, \leq\gamma$, when
$d(\xi,\xi')\leq\dt$.
\end{lemma}
\Proof By continuity it suffices to show (\ref{inner}) for
$y\in\O$. Using invariance under $G$ (and the fact that $G=G^*$),
we may assume that $y=\be$. To show that
$(\xi'|\be)\leq\gamma(\xi|\be)$, let us write
 $\xi=kP(a)\be$, for $k\in K$ and $a=a_1c_1+\ldots+a_rc_r$. Then
 the new vector $\xi''=P(a)^{-1}k^{-1}\xi'$ belongs to the fixed ball
 $B_{\dt}(\be)$. Therefore, we have
\[
(\xi'|\be)= (P(a)\xi''|\be)\leq \sqrt r\|P(a)\||\xi''|\leq
\gamma\|P(a)\|,
\]
where the last bound appears because ${\Ol{B}}_{\dt}(\be)$ is a
compact set.
 Now $P(a)$ has eigenvalues $a_i^2$ and
$a_ia_j$, and hence \Be \|P(a)\|=\max\{a^2_j,a_ia_j\} \leq
\sum_{i=1}^ra_i^2= (P(a)\be|\be)=(\xi|\be). \label{2.4}\Ee

Finally, let us remark that $(\xi|\be)$ is equivalent to $|\xi|$.
Indeed, $(\xi|\be)\leq \sqrt r |\xi| $ by Schwarz inequality.
Conversely, for $\xi=P(a)\be$, we have $$|\xi|=|\sum_{j=1}^r a_j^2
c_j |\leq \sum_{j=1}^r a_j^2= (\xi|\be)\,.$$
 \ProofEnd

\begin{lemma}
\label{g} For every $g\in G$ we have
\[
\|g\|\leq|g\be|\leq \sqrt{r}\,\|g\|.
\]
\end{lemma}
\Proof Write $g=kP(a)h$, for some $h,k\in K$ and
$a=a_1c_1+\ldots+a_rc_r$. Then, as in (\ref{2.4})
\[
\|P(a)\|\leq(a_1^4+\ldots+a_r^4)^\frac12= |P(a)\be|=|g\be|.
\]
Thus,
\[
\frac{|g\be|}{|\be|}\leq\|g\|= \|P(a)\|\leq |g\be|.
\]
\ProofEnd

\subsection{Integrals on ${\bf \Omega}$}

To conclude with this preliminary section, we list basic some
facts concerning integrals in the cone. Following \cite{FK}, we
define the {\it generalized power function in $\O$} by
\[
\Dt_\bs(x)=\Dt_1^{s_1-s_2}(x)\,\Dt_2^{s_2-s_3}(x)
\cdots\Dt_r^{s_r}(x),\quad \bs=(s_1,s_2,\ldots,s_r)\in\SC^r,\quad
x\in\O,
\]
where $\Dt_k$ are the principal minors with respect to
 a fixed Jordan frame
$\{c_1,\ldots,c_r\}$. In particular, $\Dt_\bs(x)=a_1^{s_1}\cdots
a_r^{s_r}$ when $x=a_1c_1+\ldots+a_rc_r$. The lemmas from the
previous section justify the following discretization of integrals
which we shall use often below.

\begin{proposition}
\label{discr} Let $0<\dt,R\leq R_0$ be fixed, and $\{\xi_j\}_j$ be
a $(\dt,R)$-lattice with associated Whitney decomposition
$\{E_j\}_j$. Then, for every $\bs\in\SR^r$ there exists a positive
constant $C$ such that, for any $y\in{\Ol{\O}}$ and for any
non-negative function $f$ on the cone, we have \Beas
\frac1C\,\sum_je^{-\gamma(y|\xi_j)}\Dt_\bs(\xi_j)\,\int_{E_j}
f(\xi)\,\dxi & \leq & \Ds \int_\O
f(\xi)e^{-(y|\xi)}\Dt_\bs(\xi)\,\dxi\\
 & \leq & C\,\sum_je^{-\frac1\gamma(y|\xi_j)}\Dt_\bs(\xi_j)\,\int_{E_j} f(\xi)\,\dxi\,,
\Eeas where $\gamma=\gamma(R_0,\O)$ is a constant as in
(\ref{inner}).
\end{proposition}

We shall also need the {\it gamma function in $\O$} defined from
the generalized powers. That is, given
$\bs=(s_1,s_2,\ldots,s_r)\in\SC^r$, one lets
 \Be \Ga_\O(\bs)=\int_\O
e^{-(\xi|\be)}\,\Dt_\bs(\xi)\,\dxi\,. \label{gamma} \Ee This
integral is known to converge absolutely if and only if
$\re{s_j}>(j-1)\frac{n/r-1}{r-1}$, for all $j=1,\ldots,r$.
Moreover, in such case
  \Be \Ga_\O(\bs)=(2\pi)^{\frac{n-r}2}\,
\prod_{j=1}^r\Ga(s_j-(j-1){\Ts\frac{n/r-1}{r-1}}), \label{formula}
\Ee where $\Ga$ is the classical gamma function in $\SR_+$
\cite[Ch. VII]{FK}. We shall denote $\Gamma_\O(\bs)=\Gamma_\O(s)$
when $\bs=(s,\ldots,s)$. The next formula defines the Laplace
transform of a generalized power, and can be found in \cite[p.
124]{FK}.

\begin{lemma}
\label{int} For $y\in\O$ and $\bs=(s_1,s_2,\ldots,s_r)\in\SC^r$
with $\re s_j>(j-1)\frac{n/r-1}{r-1}$, $j=1,\ldots,r$, then
\[
\int_\O e^{-(\xi|y)}\,\Dt_\bs(\xi)\,\dxi=
\Ga_\O(\bs)\,\Dt_\bs(y^{-1}).
\]
\end{lemma}

\begin{Remark}
{\rm We will sometimes write the above quantity $\Dt_\bs(y^{-1})$
in terms of the {\it rotated Jordan frame} $\{c_r,\ldots,c_1\}$.
That is, if we denote by $\Dt^*_j$, $j=1,\ldots,r$, the principal
minors with respect to this new frame, then
\[
\Dt_\bs(y^{-1})\,=\,\left[\Dt^*_{\bs^*}(y)\right]^{-1}, \quad
\forall\,\,\bs=(s_1,\ldots,s_r)\in\SC^r,
\]
 where we have set  $\bs^*:=(s_r,\ldots,s_1)$ (see \cite[p. 127]{FK}).}
\end{Remark}

Our last lemmas have to do with global and local integrability of
generalized powers. The first one is a simple consequence of our
last result and the Plancherel formula (see also \cite{BBPR}).

\begin{lemma}
\label{ial} Let $\al\in\SR$, and define
\[
I_\al(y)=\int_\SRn|\Dt(x+iy)|^{-\al}\,dx,\quad y\in\O.
\]
Then, $I_\al$ is finite if and only if $\al>\frac{2n}r-1$. In this
case, $I_\al(y)=c(\al)\,\Dt(y)^{-\al+\fnr}$.
\end{lemma}

We next establish the critical
 index for local integrability at the origin.

\begin{lemma}
\label{log} Let $\al\in\SR$ and $g_\al(\xi)=\frac{e^{-(\xi|\be)}}
{\Dt(\xi)(1+|\log\Dt_{(0,\ldots,0,1)}(\xi)|)^\al}$.
 Then, $g_\al$ is integrable if and only if $\al>1$.
\end{lemma}
\Proof This is a simple exercise using Gaussian coordinates (see
Chapter VI of \cite{FK}). Indeed, with the notation in \cite{FK},
the integral of $g_\al$ is equal to $$  c_r\int_{(0,\infty)^r}
\frac{e^{-\sum u_j^2}}
{(1+2|\log{u_r}|)^\al}\,\left[\,\prod_{j=1}^ru_j^{(r-j)d-1}\,\right]
\,du_1\ldots du_r=c'_r \int_0^{\infty} \frac{e^{-u_r^2}}
{(1+2|\log{u_r}|)^\al}\,\frac{du_r}{u_r}. $$

\ProofEnd

Finally, we conclude with the critical index for integrability at
infinity.

\begin{lemma}
\label{log2} Let $\al,\dt\inR$, $\beta>-1$ and
$g_{\al,\beta,\dt}(y)=\frac{\Dt^\beta(y)}
{\Dt^\al(y+\be)(1+\log\Dt(y+\be))^\dt}$.
 Then, $g_{\al,\beta,\dt}$ is integrable if and only if $\al-\beta>\frac{2n}r-1$
 or  $\al-\beta=\frac{2n}r-1$ and $\dt>1$.
\end{lemma}
\Proof This time we use the ``polar coordinates'' of the cone
\[
y=k(e^{t_1}c_1+\ldots+e^{t_r}c_r),\quad t_1<t_2<\ldots<t_r,\quad
k\in K
\]
(see \cite[pag. 105]{FK}). Then,
$\Dt(y+\be)=\prod_{j=1}^r(e^{t_j}+1)$, and
\[
\int_\O g_{\al,\beta,\dt}(y)\,dy= c\,\int_{-\infty}^\infty
\int_{-\infty}^{t_r}\cdots\int_{-\infty}^{t_2}
\frac{e^{(t_1+\cdots+t_r)(\fnr+\beta)}\,\prod_{j<k}\left(\sh({\Ts\frac{t_k-t_j}2})\right)^d}
{\prod_{j=1}^r(e^{t_j}+1)^\al\,(1+\sum_{j=1}^r\log(1+e^{t_j}))^\dt}\,dt_1\ldots
dt_r,
\]
where $d=\dim V_{j,k}=2(\fnr-1)/(r-1)$. For the necessary
condition, we can consider only the case when
$\al-\beta=\frac{2n}r-1$ and $\dt=1$. Moreover, we restrict the
region of integration so that $\sh t\geq ce^t$, and obtain \Beas I
& \geq & c\, \int_{2^r}^\infty
\int_{2^{r-1}}^{2^{r-1}+1}\cdots\int_2^3
\frac{e^{(t_1+\cdots+t_r)(1-\fnr)}}{1+t_r}\, \prod_{j<k} e^{\frac
d2(t_k-t_j)}\, dt_1\ldots dt_r\\
 & \geq & c'\,\int_{2^r}^\infty
\frac{e^{t_r(1-\fnr)}}{1+t_r}\,e^{\frac{d(r-1)}2t_r}\,
dt_r=\infty. \Eeas

To estimate from above,
 we  use the bound
\[
\prod_{j<k}\sh({\Ts\frac{t_k-t_j}2})\leq
\prod_{j<k}e^{\frac{(t_k-t_j)}2}=\prod_{j=1}^r\,e^{(j-1-\frac{r-1}2)t_j}.
\]
Then the integral $I$ is bounded by the product
\[
\prod_{j=1}^{r-1}\int_{-\infty}^{+\infty}\frac{e^{(d(j-1)+\beta+1)t_j}}{(1+e^{t_j})^\alpha}dt_j
\times
\int_{-\infty}^{+\infty}\frac{e^{(d(r-1)+\beta+1)t_r}}{(1+e^{t_r})^\alpha(1+\log
(1+e^{t_r}))^\delta }dt_r\,.\] Each integral is convergent at
 $-\infty$  since
$\beta+1>0$. We use the conditions on $\alpha$, $\beta$ and
$\delta$ to conclude easily for the integrability at $+\infty$.
 \ProofEnd

\section{Besov spaces with spectrum in $\O$}
\setcounter{equation}{0} \setcounter{footnote}{0}
\setcounter{figure}{0}

\subsection{The Littlewood-Paley decomposition}

Through the rest of the paper, $\{\xi_j\}$ will be a fixed
$(\dt,R)$-lattice in $\O$ with $\dt=\frac12$ and $R=2$. We can
easily construct a smooth partition of the unity associated with
the covering $B_j=B_1(\xi_j)$. For this, we choose a real function
 $\phi_0\in C^\infty_c(B_2(\be))$ such that
\[
0\leq \phi_0 \leq 1, \mand \phi_0|_{B_1(\be)}\equiv1.
\]
We  write each point $\xi_j=g_j\be$, for some fixed $g_j\in G$
(which, for simplicity,
 we take self-adjoint).
Then, we can define $\phi_j(\xi):=\phi_0(g_j^{-1}\xi)$, so that
 \Be
 \phi_j\in C^\infty_c(B_2(\xi_j)),\quad
0\leq \phi_j\leq 1
 \mand
 \phi_j|_{B_j}\equiv1.
 \label{phij}
 \Ee
We assume that $\xi_0=\be$, so that there is no ambiguity of
notations. By the finite intersection property, there exists a
constant $c>0$ such that
\[
{\Ts\frac1c}\leq \Phi(\xi):=\sum_j\phi_j(\xi) \leq c.
\]

\begin{proposition}
\label{unity} In the conditions above, let $\hpsi_j=\phi_j/\Phi$.
Then \Benu \item $\hpsi_j\in C^\infty_c(B_2(\xi_j))$;

\item $0\leq \hpsi_j\leq 1,\mand \sum_j\hpsi_j(\xi)=1,
\,\,\forall\,\xi\in\O$;

\item $\psi_j$ are uniformly bounded in $L^1(\SRn)$; in
particular, there exists a constant $C>0$ such that
 \Be
 \|f*\psi_j\|_p\leq C\,\|f\|_p,\quad
 \forall\,f\in L^p(\SRn),\,\,\,\forall\,j,
 \,\,\,1\leq p\leq\infty.
 \label{mult1}
 \Ee
\Eenu
\end{proposition}
\Proof The first two statements are clear. For the last one, note
first
 that
\[
\|\psi_j\|_{L^1}= \|\cF^{-1}(\phi_0(g_j^{-1}\,\cdot)/\Phi)\|_{L^1}
=\|\cF^{-1}(\phi_0/\Phi(g_j\,\cdot))\|_{L^1}.
\]
Now, when $\xi\in B_2(\be)$ we can write
\[
\Phi(g_j\xi)=\sum_{k\in J_j}\phi_0(g_k^{-1}g_j\xi),
\]
where $J_j=\{k\mid B_2(\xi_k)\cap B_2(\xi_j)
 \not=\emptyset\}$ is
a finite set with at most $N$ elements by the finite intersection
property.
 Further, we claim that the following uniform estimate holds
 true:
 \Be
\|g^{-1}_kg_j\|\leq C, \quad\mbox{when}\,\,\, k\in
J_j,\,\,\forall\,j. \label{unif} \Ee Indeed, since
$d(g_k^{-1}g_j\be,\be)=d(\xi_j,\xi_k)\leq 4$, by Lemmas \ref{g}
and \ref{l6},
    \[ \|g^{-1}_kg_j\|  \sim
|g^{-1}_kg_j\be| \sim|\be|\sim1.
\]

From (\ref{unif}) the proposition  follows easily. Indeed,
integrating by parts we have \Bea
\cF^{-1}(\phi_0/\Phi(g_j\,\cdot))(x) & = &
\int_{B_2(\be)}\,e^{i(x|\xi)}\,\frac{\phi_0(\xi)}
{\Phi(g_j\xi)}\,d\xi\nonumber\\ & = &
\int_{B_2(\be)}\,e^{i(x|\xi)}\,
\frac{D^L(\phi_0/\Phi(g_j\cdot))(\xi)}{(-|x|^2)^L}
\,d\xi,\label{parts} \Eea where $D^L$ denotes a power of the
Laplacian. All functions $D^L(\phi_0/\Phi(g_j\cdot))$ are bounded.
Thus, choosing $L=0$ for $|x|\leq1$, and $L>\frac{n}2$ for
$|x|>1$, we can majorize $\cF^{-1}(\phi_0/\Phi(g_j\,\cdot))$
uniformly in $j$ by an integrable function, and this establishes
the result.\ProofEnd
\smallskip

In this paper we shall mainly be concerned with  Besov-type
seminorms derived from the couple $\{\xi_j,\psi_j\}$ as in
(\ref{bpqnorm}). That is, for $\nu\in\SR$, $1\leq p,q\leq\infty$,
and $f\in\cS'(\SR^n)$ we let
 \Be \|f\|_\Bpq:=\left\{\Ba{lcl}
 \left[\sum_j\Dt^{-\nu}(\xi_j)\,\|f*\psi_j\|_p^q\,\right]^\frac1q& & \mbox{if }q<\infty\\
 & & \\
 \sup_j\Dt^{-\nu}(\xi_j)\,\|f*\psi_j\|_p& & \mbox{if }q=\infty\,.
 \Ea\right. \label{bpqnorm1}\Ee
 We shall make use of the fact that these
seminorms do not actually depend on the choice of the lattice
$\{\xi_j\}$ or the test functions $\psi_j$. Moreover, they can as
well be defined with test functions which are not normalized as in
the previous proposition. That is, we may replace $\{\psi_j\}$ by
any family
    \Be
    \hchi_j(\xi):=\hchi(g_j^{-1}\xi),
    \label{chi}
    \Ee
defined from an arbitrary $\hchi\in C^\infty_c(B_4(\be))$ so that
$0\leq\hchi\leq1$ and  $\hchi$ is identically $1$ in $B_2(\be)$.
These and other elementary equivalences are stated and proved in
the following lemma.

\begin{lemma}
\label{besov-norm-lem} Let $\{\xi_j,\psi_j\}$ be as at the
beginning of this section, and fix $\nu\in\SR$ and $1\leq
p,q\leq\infty$. Then, for any other $(\dt,R)$-lattice $\{\txi_j\}$
with associated Littlewood-Paley functions $\{\tpsi_j\}$, and for
any family $\{\chi_j\}$ as in (\ref{chi}), we have the
equivalences
$$\Bigl[\sum_j\Dt^{-\nu}(\xi_j)\,\|f*\psi_j\|_p^q\,\Bigr]^\frac1q\,\sim\,
\Bigl[\sum_j\Dt^{-\nu}(\txi_j)\,\|f*\tpsi_j\|_p^q\,\Bigr]^\frac1q\,,$$
$$\Bigl[\sum_j\Dt^{-\nu}(\xi_j)\,\|f*\psi_j\|_p^q\,\Bigr]^\frac1q\,\sim\,
\Bigl[\sum_j\Dt^{-\nu}(\xi_j)\,\|f*\chi_j\|_p^q\,\Bigr]^\frac1q\,,$$
for all $f\in\cS'(\SR^n)$. Moreover, when $g\in G$ and $q<\infty$,
it holds the equivalence
$$\sum_j\Dt^{-\nu}(\xi_j)\,\|(f\circ g)*\psi_j\|_p^q\,\sim\, \Dt
(g\be)^{-\fnr\frac qp-\nu}
\sum_j\Dt^{-\nu}(\xi_j)\,\|f*\psi_j\|_p^q.$$
\end{lemma}

\Proof We just consider the case $q<\infty$, the modifications for
$q=\infty$ being trivial. For the first part, we can write, for
each $j$,
    \[
    \hpsi_j=\sum_k\hpsi_j\htpsi_k,
    \]
where the index $k$ runs through a set $J_j=\{k\mid
B_{R\dt}(\txi_k)\cap B_2(\xi_j)\not=\emptyset\}$ of at most
$N=N(\dt,R,\O)$ elements. Then, using (\ref{mult1}) and Lemma
\ref{l3} we have
    \Beas
    \sum_j\Dt^{-\nu}(\xi_j)\,\|f*\psi_j\|_p^q & \leq &
C\,\sum_j\sum_{k\in J_j}\Dt^{-\nu}(\xi_j)\,\|f*\tpsi_k\|_p^q\\
    & \leq &
C'\,\sum_k\Dt^{-\nu}(\txi_k)\,\|f*\tpsi_k\|_p^q.
    \Eeas
The converse inequality follows similarly. For the second
equivalence in the lemma, the fact that $\hchi_j\hpsi_j=\hpsi_j$
implies immediately the left inequality. A similar use of the
finite intersection property as we did above gives the right hand
side.

Finally, for our last statement, it is sufficient to prove an
inequality of the form ``$\lesssim$'', the converse inequality
``$\gtrsim$'' following after replacing $g$ by its inverse. Now,
using a first change of variables, and the fact that the
determinant of the transformation $g$ in $\SR^n$ is equal to
$\Dt(g\be)^{\frac nr}$, we are linked to consider the $L^p$-norms
 of the functions
 $$\Dt(g\be)^{-(1+\frac 1p)\frac nr}\,f*(\psi_j\circ g^{-1})
 =\Dt(g\be)^{-(1+\frac 1p)\fnr}\sum _k f*\psi_k*(\psi_j\circ g^{-1})\,. $$
For each fixed $j$ this last sum has at most $N$ terms, since the
Fourier transform of $\psi_k*(\psi_j\circ g^{-1})$ is non zero
only if $d(\xi_k,g^*\xi_j)<4$. So $$\|f*(\psi_j\circ
g^{-1})\|_p\leq C \Dt(g\be)^{\frac nr}\sum_{k ;
d(g^*\xi_k,\xi_j)<4}\|f*\psi_k\|_p,$$ the factor $\Dt(g\be)^{\frac
nr}$ appearing as the determinant of the transformation $g$ in the
computation of the $L^1$ norm of $\psi_j\circ g^{-1}$. Now, when
$d(g^*\xi_k,\xi_j)<4$, then $\Dt(\xi_j)$ is equivalent to
$\Dt(g^*\xi_k)=\Dt(g\be) \Dt(\xi_k)$. Thus, we conclude $$
\sum_j\Dt^{-\nu}(\xi_j)\,\|f*(\psi_j\circ g^{-1})\|_p^q\leq C
\Dt(g\be)^{\frac {nq}r-\nu}\sum_j\sum_{k ;
d(g^*\xi_k,\xi_j)<4}\Dt^{-\nu}(\xi_k)\|f*\psi_k\|_p^q.$$ We get
the required inequality multiplying by $\Dt(g\be)^{-q(1+\frac
1p)\frac nr}$ and summing first in the $j$ indices.
 \ProofEnd
\smallskip

 Recall
now that $\So$ denotes the space of Schwartz functions $f$ on
$\SRn$ with $\supp\hf\subset\clO$. The next proposition gives the
Littlewood-Paley decomposition $f=\sum_j f*\psi_j$ for functions
in $f\in\So$, and relates it with the Besov space norm.

\begin{proposition}
\label{So} Every $f\in\So$ admits a Littlewood-Paley decomposition
$f=\sum_j f*\psi_j$ with convergence in $\cS(\SRn)$. Further, for
every $\nu\inR$, $1\leq p,q\leq\infty$ there is a constant
$C=C(p,q,\nu)>0$ and an integer $\ell=\ell(p,q,\nu)\geq0$ so that
\Be \|f\|_{\Bpq}=\left[\sum_j\Dt^{-\nu}(\xi_j)\,\|f*\psi_j\|_p^q
\right]^\frac1q \leq C\,p_{\ell}(\hf)<\infty,
\quad\forall\;f\in\So, \label{So1} \Ee where
$p_{\ell}(\phi)=\sup_{|\al|\leq \ell}\sup_{\xi\in\SRn}
(1+|\xi|)^\ell\,|\partial^\al\phi(\xi)|$ denotes a Schwartz
seminorm.
\end{proposition}

The proof depends on a lemma which gives appropriate estimates for
test functions in $\So$.

\begin{lemma}
\label{schwartz} Let $N,M\geq 0$. Then, there is a constant
$C=C(N,M)>0$ and an integer $\ell=\ell(N,M)\geq0$ such that for
every $f\in\So$ \Benu \item
    $|\hf(\xi)|\leq\,C\,p_\ell(\hf)\,\Dt^M(\xi)\,(1+|\xi|)^{-N},
    \quad\forall\,\xi\in\O;$
 \item
    If $1\leq p\leq \infty$,
    $\|f*\psi_j\|_p\leq\,C\,p_\ell(\hf)\,\Dt(\xi_j)^{M+\fnr\frac1{p'}}
    \,(1+|\xi_j|)^{-N},
    \quad\forall\,j.$
    \Eenu
\end{lemma}
\Proofof{Lemma \ref{schwartz}}

For the first statement,  it suffices to show that for every
$f\in\So$ and $M\geq1$ there is $M'\geq1$ so that
    \Be
    |\hf(\xi)|\leq
    \,C\,p_{M'}(\hf)\,\Dt^M(\xi),\quad\mbox{whenever}\quad\Dt(\xi)\leq1,
    \quad\xi\in\O.
    \label{sch1}
    \Ee
Indeed, we then write it for $D^Nf$ to get the full statement. So,
let us prove (\ref {sch1}).  Let $\xi\in\O$ be fixed, and choose
$\xi_0\in\bO$ so that dist$(\xi,\bO)=|\xi-\xi_0|$. Since
$\supp\hf\subset\clO$, we  have $\partial^\al\hf(\xi_0)=0$, for
every multi-index $\al$. Thus, given $M\geq1$ there is a constant
$C=C(M)$ such that $|\hf(\xi)|\leq C\,p_M(\hf)\,|\xi-\xi_0|^M$. We
claim that
 $|\xi-\xi_0|\leq\Dt(\xi)^\frac1r$, which
  will clearly establish (\ref{sch1}).

To show our claim, we may assume that $\xi=P(a)\be$, where
$a=a_1c_1+\ldots+a_rc_r$. Suppose also that
$a_1=\min\{a_1,\ldots,a_r\}$. Then
  \Beas
    |\xi-\xi_0| & = & \mbox{dist }(\xi,\bO)\,\leq\,
    |\xi-(a_2^2c_2+\ldots+a_r^2c_r)|\\
     & = &
a_1^2\,\leq\,(a_1^2\,\cdots\,a_r^2)^\frac1r=\Dt(\xi)^\frac1r.
     \Eeas

Let us now prove the second statement in Lemma \ref{schwartz}. It
is sufficient to prove the same inequality, with the system
$\chi_j$ instead of $\psi_j$.  Given $f\in\So$, we proceed as in
(\ref{parts}):
 \Bea
    f*\chi_j(x) & = & \int_\O e^{i(x|\xi)}\,
    \hf(\xi)\hchi(g_j^{-1}\xi)\,d\xi\nonumber\\
    & = & \Dt^\fnr(\xi_j)\,\int_\O e^{i(g_jx|\xi)}\,
    \hf(g_j\xi)\hchi(\xi)\,d\xi\nonumber\\
    & = & \Dt^\fnr(\xi_j)\,\int_{B_2(\be)} e^{i(g_jx|\xi)}\,
   {{D^{L}(\hf(g_j\xi)\hchi(\xi))}\over{(-|g_jx|^2)^L}}
    \,d\xi.\label{p17}
    \Eea
The estimates in the first part, together with Lemmas \ref{l3},
\ref{l6} and \ref{g}, imply that, on the invariant ball
$B_2(\be)$, \Beas |D^{L}(\hf(g_j\xi))| & \leq & C\,
(1+\|g_j\|)^{2L}\,\sum_{|\al|\leq 2L}\,|(\partial^\al\hf)(g_j\xi)|
\\ & \leq & C'\,p_\ell(\hf)\,
\frac{\Dt^M(\xi_j)}{(1+|\xi_j|)^N}, \Eeas for some integer
$\ell=\ell(M,N,L)$. Therefore
  \[
   |f*\chi_j(x)|\leq\,C\,p_\ell(\hf)\,\Dt^\fnr(\xi_j)
\,\frac{\Dt^M(\xi_j)}{(1+|\xi_j|)^N}\,
\frac1{(1+|g_jx|^2)^L},\quad x\in\SRn.
    \]
Taking $L^p$-norms and changing variables, we conclude with \Be
\|f*\chi_j\|_p\,\leq\,C\,p_\ell(\hf)\, \frac{\Dt(\xi_j)^{M+\frac
nr\frac1{p'}}}{(1+|\xi_j|)^N}.
    \label{mnp}
\Ee \ProofEnd

\Proofof{Proposition \ref{So}} Once we show the convergence of the
series $\sum_j f*\psi_j$, the fact that the sum equals $f$ is
immediate from $\supp\hf\subset\clO$ and $\sum_j\hpsi_j=\chi_\O$.
Now, from the previous lemma and Proposition \ref{discr}, we have
\Beas \sum_j\|\hf\hpsi_j\|_\infty & = & \sum_j\|\hf(g_j\cdot)
\,\frac{\phi_0}{\Phi(g_j\cdot)}\|_\infty\\
 & \leq & C_f\,
\sum_j\Dt^M(\xi_j)\,(1+|\xi_j|)^{-N}\leq \,C'_f\,\int_\O
    \frac{\Dt^M(\xi)}{(1+|\xi|)^N}\,\dxi,\Eeas
where the last integral is finite for $N,M$ large enough. A
similar argument applies to
$\|(1+|\xi|)^L\partial^\al(\hf\hpsi_j)\|_\infty$, establishing our
claim.

For the second assertion in the proposition, we use the second
estimate in Lemma \ref{schwartz}. Assuming $q<\infty$ (otherwise
the estimate is trivial) we have
 \Bea
    \|f\|^q_\Bpq & \leq & \,C\,p_\ell(\hf)^q\,\sum_j\Dt^{-\nu}(\xi_j)\,
    \frac{\Dt(\xi_j)^{Mq+\frac{n}r\frac{q}{p'}}}
{(1+|\xi_j|)^{Nq}}\nonumber\\
     & \leq & \,C'\,p_\ell(\hf)^q
    \int_\O
    \frac{\Dt(\xi)^{Mq+\frac nr\frac{q}{p'}-\nu}}{(1+|\xi|)^{Nq}}\,\dxi,
\label{doub} \Eea which  is finite for a sufficiently large choice
of $N,M$. \ProofEnd

We observe that we have strongly used the assumption on the
support of $\hat f$.
 The next proposition gives the sharp region of indices
$\nu,p,q$ for which general Schwartz functions have finite
$\Bpq$-seminorms.
  We state this
fact separately since such conditions will appear in the sequel in
relation with the index $\tqn$.

\begin{proposition}
\label{Rsch} Let $\nu\inR$, $1\leq p,q\leq\infty$ be such that \Be
\frac q{p'}\,\fnr>\nu+\fnr-1\label{vpq}\Ee (or
$\frac1{p'}\fnr\geq\nu$ when $q=\infty$). Then, there exist
$C=C(p,q,\nu)>0$ and an integer $\ell=\ell(p,q,\nu)\geq0$ so that:
\Be \|f\|_{\Bpq} \leq C\,p_{\ell}(\hf)<\infty,
\quad\forall\;f\in\cS(\SR^n). \label{So2} \Ee Moreover, this
property can hold for all $f\in\cS(\SR^n)$ only if $\nu,p,q$
satisfy (\ref{vpq}).
\end{proposition}
\Proof As before we assume $q<\infty$, with obvious modifications
when $q=\infty$. First one observes that, when $f\in\cS(\SRn)$,
the conclusion of Lemma \ref{schwartz} is still valid with $M=0$.
Thus, a similar reasoning as in (\ref{doub})
 gives
\[
\|f\|_\Bpq\leq \,C \,p_\ell(\hf)\, \left[\,\int_\O
\frac{\Dt^{\frac
q{p'}\fnr-\nu}(\xi)}{(1+|\xi|)^{Nq}}\,\dxi\,\right]^\frac1q, \quad
f\in\cS(\SRn),
\]
and this integral is  finite under the condition $\frac
q{p'}\fnr>\nu+\fnr-1$. To show that this condition is critical,
take any $f\in\cS(\SR^n)$ such that $\hf$ is identically $1$ in
the Euclidean ball centered at $0$ of radius $1$. Then, for such
an $f$, one has the bound from below
$$\|f\|_\Bpq^q \geq c \sum_{j;
|\xi_j|<c}\Dt(\xi_j)^{-\nu}\|\chi_j\|_p^q \geq c''\int_{\xi\in\O;
|\xi|<c'}\Dt^{\frac q{p'}\fnr-\nu}(\xi)\dxi,$$ and this  last
integral is infinity unless $\frac q{p'}\fnr>\nu+\fnr-1$.
 \ProofEnd

\subsection{Properties of Besov spaces}

Given a closed set $F\subset\SRn$, we shall denote by
$S'_F=S'_F(\SRn)$ the space of tempered distributions with Fourier
transform supported in $F$.  Recall the expression of the
``seminorm'' $\|f\|_\Bpq$ in (\ref{bpqnorm}), and observe that a
distribution $f\in\Sop$ satisfies $\|f\|_\Bpq=0$ if and only if
$f\in\cS'_{\bO}$. This leads to the following natural definition
of Besov spaces with spectrum in $\O$:

\begin{definition}
Given $\nu\inR$, $1\leq p,q<\infty$, we define $\Bpq$ as the space
of equivalence classes of tempered distributions
\[
\Bpq=\{f\in\Sop\mid
\|f\|_\Bpq<\infty\}\,\left/\,\cS'_{\bO}\right..
\]
\end{definition}

It follows right away from Lemma \ref{besov-norm-lem} that $\Bpq$
does not depend on the choice of $\{\psi_j\}$ or the lattice
$\{\xi_j\}$. Moreover, $\Bpq$ is invariant under the action of
$G$, and \Be \label {g-action}
 \|f\circ g\|_\Bpq\sim \Dt(g\be)^{-\frac n{rp}- \frac{\nu}q}
 \|f\|_\Bpq.\Ee
Before collecting in the next proposition other basic properties
of the spaces $\Bpq$, we make a minor comment on the notation used
below.

\begin{notation}\label{notation}
{\rm Throughout this paper, the standard action of tempered
distributions over Schwartz functions will be denoted by
\[
(f,\phi)=\int f\,\phi,\quad f\in\cS'(\SR^n),\;\phi\in\cS(\SR^n).
\]
 For convenience,  we shall often use the \emph{anti-linear pairing}:
 \[
\langle f,\phi\rangle:=(f,\bar\phi)=\int_{\SR^n} f\,{\Ol{\phi}},
\quad f\in\cS'(\SR^n),\;\phi\in\cS(\SR^n).
\]
This has the notational advantage of a simple Plancherel identity:
$\langle f,\phi\rangle=\langle \hf,\hphi\rangle$, leading to a
natural pairing between $f\in\Sop(\SR^n)$ and $\phi\in\So$ (rather
than using $(f,\phi)=(\hf,\hphi(-\,\cdot\,))$, which requires to
deal with $\phi\in\cS_{-\O}$).
 }\end{notation}

With the above considerations we have the following lemma.

\begin{lemma}
\label{obser} Let $\nu\inR$, $1\leq p\leq\infty$ and $1<q<\infty$.
Then, there exists $\ell=\ell(\nu,p,q)\geq0$ so that, for every
distribution $f\in\cS'(\SRn)$ with $\|f\|_\Bpq<\infty$, we have
\Be |\langle
f,\phi\rangle|\leq\,C\,\|f\|_\Bpq\,\|\phi\|_{B_{-\nu{q'}/q}^{p',q'}}
\leq C'\,\|f\|_\Bpq\,p_\ell(\hphi),\quad
\forall\;\phi\in\So.\label{bholder} \Ee When $q=1$ or $q=\infty$,
the same holds replacing $\|\phi\|_{B_{-\nu{q'}/q}^{p',q'}}$ by
$\|\phi\|_{B_{-\nu}^{p',q'}}$.
\end{lemma}
 \Proof

Remember that $\hpsi_j=\hpsi_j\hchi_j$ for all $j$. Therefore,
using Proposition \ref{So} and Lemma \ref{besov-norm-lem}
 \Be \vert \langle f,\phi\rangle\vert  \leq
\sum_j\vert\langle f*\psi_j,\phi*\chi_j\rangle\vert \, \leq \,
\sum_j\|f*\psi_j\|_p\|\phi*\chi_j\|_{p'} \leq
C\,\|f\|_\Bpq\,\|\phi\|_{B_{-\nu{q'}/q}^{p',q'}}. \Ee Finally,
observe that
$\|\phi\|_{B_{-\nu{q'}/q}^{p',q'}}=\|\phi\|_{B_{\nu(1-q')}^{p',q'}}$,
which, for $\phi\in\So$, is bounded by a Schwartz seminorm  by
Proposition \ref{So}. The modifications for $q=1,\infty$ are
obvious. \ProofEnd

\begin{proposition}
\label{besovprop} Let $1\leq p,q<\infty$ and $\nu\in\SR$. Then
\Benu \item $\Bpq$ is a Banach space. \item The space
$\Do:=\{f\in\cS(\SRn)\mid \supp\hf$ is
compact in $\O\}$ is dense in $\Bpq$. Moreover, for every class
$f+\Sbop$ in $\Bpq$, the series $\sum_jf*\psi_j$ converges to (the
class of) $f$ in the space $\Bpq$. \Eenu
\end{proposition}
\Proof

Suppose $\{f_m\}_m$ is a Cauchy sequence of distributions in
$\Sop$ for the $\Bpq$-seminorm. Then, from the previous lemma it
follows that $\langle f_m,\phi\rangle$ converges for every
$\phi\in\So$, and moreover, it defines a continuous (anti-linear)
functional in $\So$. We can extend it to $\So\oplus\cS_{\O^c}$ by
letting it be identically zero in the second summand, and finally
extend it to the whole Schwartz space $\cS(\SRn)$ by the
Hahn-Banach theorem. This gives a tempered distribution $f\in\Sop$
which in particular satisfies
\[
f*\psi_j(x)=\lim_{m\to\infty}
f_m*\psi_j(x),\quad\forall\;x\in\SRn,\quad \forall\;j.
\]
Therefore, by Fatou's lemma
\[
\|f\|_\Bpq\leq \liminf_{m\to\infty}\|f_m\|_\Bpq<\infty,
\]
and in a similar fashion $\lim_{m\to\infty}\|f-f_m\|_\Bpq=0$. This
shows that $\Bpq$ is a Banach space.

For the density, let $f$ be a fixed distribution in $\Sop$ with
$\|f\|_\Bpq<\infty$. We shall show that $f$ is the $\Bpq$-limit of
the partial sums of the series $\sum_{j}f*\psi_j$. Remark that
each finite sum belongs to $L^p(\SRn)$, and therefore  can be
approached by a Schwartz function with Fourier transform supported
in a compact set of $\O$, justifying the density of $\Do$. Now, it
is easily seen that partial sums (for any order) constitute a
Cauchy sequence in $\Bpq$. Since $\Bpq$ is a Banach space, they
converge to a distribution $u\in\Sop$. It remains to show that $u$
and $f$ belong to the same equivalence class in $\Sop/\Sbop$,
which is an immediate consequence of the fact that
$f*\psi_j=u*\psi_j$. \ProofEnd

For the duality of the spaces $\Bpq$, recall the H\"older type
inequality in (\ref{bholder}): \[ |(\bar f,g)|=|\lan fg|\leq
C\,\|f\|_\Bp\,\|g\|_\Bpq,\quad g\in\Do,
\]
valid for every $f\in\Sop$ with $\|f\|_\Bp<\infty$. Observe that
$g\in\Do\mapsto(\bar f,g)$ is a linear functional in $\Do$ which
 depends only on the equivalence class $f+\Sbop$. Thus, by the
above inequality and the density of $\Do$, it defines a continuous
linear functional $\Phi_f$ in $\Bpq$. Further, if $\Phi_f=0$, then
$(\bar f,g)=0$, $\forall$ $g\in\Do$, and necessarily $f\in\Sbop$.
Thus, the correspondence \Be
f+\Sop\in\Bp\longrightarrow\Phi_f\in(\Bpq)^* \label{dual1}\Ee
 is well-defined and
injective.

\begin{proposition}
\label{bpqdual} Let $\nu\inR$ and $1\leq p,q<\infty$. Then, the
mapping in (\ref{dual1}) is an anti-linear isomorphism of Banach
spaces.
\end{proposition}

\Proof By the previous comments it suffices to show that for every
$\Phi\in(\Bpq)^*$, there exists a distribution $f\in\Sop$ such
that \Be
\Phi(g)=(\bar{f},g),\quad\forall\;g\in\Do\quad\mbox{and}\quad
\|f\|_\Bp\leq C\,\|\Phi\|. \label{for2} \Ee Now, since
$\Do\subset\So\hookrightarrow\Bpq$, by Hahn-Banach we can extend
continuously $\Phi$ to $\cS(\SR^n)$, and find a tempered
distribution $f\in\Sop$ such that $\Phi(g)=(\bar f,g)$,
$\,\forall\,g\in\Do$.

We now claim that each $f*\psi_j$, which a priori is only a smooth
function with polynomial growth, does belong to $L^{p'}(\SRn)$,
and  moreover, the sequence of their $L^{p'}$ norms belongs to the
suitable space of sequences. Indeed,  for every finite sequence
$g_j\in\cS(\SRn)$ with $\sum_j\Dt(\xi_j)^{-\nu}\|g_j\|_p^q\leq1$
we have \Beas |\sum_j\lan{ f*\psi_j}{g_j}| & = & |\Phi(\sum_j
g_j*\psi_j)|\, \leq\, \|\Phi\|\,\|\sum_j g_j*\psi_j\|_\Bpq \leq
C\,\|\Phi\|. \Eeas The constant $C$ depends only on the number $N$
in the finite intersection property, the constant $\gamma$ related
to the variation of the function $\Dt$ inside an invariant ball of
radius $2$, and the $L^1$ norm of the $\psi_j$'s. Since the
constant is independent of the finite set of indices,
 (\ref{for2}) follows. We do not give the details of the proof,
 since it is completely analogous to the one of Lemma
 \ref{besov-norm-lem}.
\ProofEnd

 Let us remark that, for two classes of
tempered distributions $f+\Sbop$ in $\Bpq$ and $g+\Sbop$ in $\Bp$,
the duality pairing can also be expressed as
 \Be
\Phi_f(g+\Sop)=
\sum_j\lan{f*\psi_j}{g*\chi_j},\label{crochet}
\Ee where the series converges absolutely by (\ref{bholder}). This
representation is sometimes convenient, and of course, independent
on the choice of $\{\psi_j,\chi_j\}$.

\subsection{The $\Box$ operator and Besov multipliers}
\smallskip

Next we describe some analytic properties of the spaces $\Bpq$.
The first one concerns the role of the generalized wave operator
$\Box$ (introduced in (\ref{box})) as a natural isomorphism
between these spaces.
 Below we shall be interested in fractional and
 negative powers of $\Box$, which can be defined
by the rule \Be \Box^\beta f=\cF^{-1} (\Dt^\beta\hf), \label{Ib}
\Ee at least for distributions $f\in\cS'(\O)$ so that  $\supp\hf$
is compact in $\O$. Our next result is a more general version than
(2) in Theorem \ref{thb}.

\begin{proposition}
\label{ib} Let $\nu,\beta\inR$ and $1\leq p,q<\infty$. Then there
is a constant $C>0$ such that for every distribution
$f\in\cS'(\SR^n)$ with $\supp\hf$ compact in $\O$ \Be
\frac1C\,\|f\|_\Bpq\leq\|\Box^\beta f\|_{B^{p,q}_{\nu+\beta
q}}\leq C\|f\|_\Bpq. \label{boxb} \Ee In particular, $\Box^\beta$
extends to an isomorphism $\tBox^\beta\colon\Bpq\to
B^{p,q}_{\nu+\beta q}$.
\end{proposition}

\Proof Indeed, given $f\in\cS'(\SRn)$ with $\supp\hf$ compact in
$\O$ we have
\[
\|\Box^\beta f\|_{B^{p,q}_{\nu+\beta q}}^q=
\sum_j\Dt(\xi_j)^{-(\nu+\beta q)} \,\|\cF^{-1}(\hf
\hpsi_j\Dt^\beta)\|^q_p \,\leq\Tquad
\]
\[
\leq\,\sum_j\Dt(\xi_j)^{-(\nu+\beta q)}\|f*\psi_j\|^q_p\,
\|\cF^{-1}(\hchi_j\Dt^\beta)\|^q_1.
\]
Using $\Dt(g\xi)=\Dt(g\be)\Dt(\xi)$, for $g\in G$, we have
\[
\|\cF^{-1}(\hchi_j\Dt^\beta)\|_1=\Dt^\beta(\xi_j)\,
\|\cF^{-1}(\hchi\Dt^\beta)\|_1=c_\beta\Dt^\beta(\xi_j),
\]
from which (\ref{boxb}) follows easily. To extend $\Box^\beta$ to
the space $\Bpq$ one proceeds by density. More precisely, given
any $f\in\Sop$ with $\|f\|_\Bpq<\infty$, we denote by $\tBox^\beta
f$ a representative from the equivalence class of
$\sum_j\Box^\beta(f*\psi_j)$, which by (\ref{boxb}) (and
Proposition \ref{besovprop}) is a Cauchy series in the
$B^{p,q}_{\nu+\beta q}$-seminorm. Observe that
$\tBox^\beta\Sbop=0$ (or its equivalence class), while by
uniqueness of the extension, $\tBox^\beta$ does not depend on the
Littlewood-Paley functions $\{\psi_j\}$.\ProofEnd

\smallskip

A further step in the previous idea leads to a \emph{functional
calculus in $\Bpq$} based on the operator $\Box$. Let $m\in
C^\infty(0,\infty)$ be a Mihlin-type multiplier in 1 dimension.
That is, there is a constant $C=C(m)$ so that \Be
\sup_{\xi>0}|\xi|^k\,|m^{(k)}(\xi)|\leq C,\quad \forall\,
k=0,1,\ldots \label{mih} \Ee Then, it makes sense to define the
operator $m(\Box)$ by
\[
m(\Box)(f)=\cF^{-1}(m(\Dt)\hf),
\]
at least for $f\in\cS'(\SRn)$ with $\supp\hf$ compact in $\O$.
Observe that a typical example is given by the imaginary powers
$m(\xi)=\xi^{i\gamma}$ for $\gamma\inR$. Then we have the
following:

\begin{proposition}
Let $\nu\inR$, $1\leq p,q<\infty$ and $m\in C^\infty(0,\infty)$
satisfying (\ref{mih}). Then, there is a constant $C=C(m)$ such
that
\[
\|m(\Box)(f)\|_\Bpq\leq C_m\,\|f\|_\Bpq,\quad f\in\Do.
\]
In particular, $m(\Box)$ extends to a bounded operator in $\Bpq$.
\end{proposition}

\Proof For the proof it suffices to estimate \Be
\|m(\Box)(f)*\psi_j\|_p\leq \|f*\psi_j\|_p\,
\|\cF^{-1}(m(\Dt)\hchi_j)\|_1,\quad\mbox{for}\quad f\in\Do.
\label{jus1} \Ee Now, \Bea \|\cF^{-1}(m(\Dt)\hchi_j)\|_1 & = &
\|\cF^{-1}(m(\Dt(\xi_j)\Dt)\hchi)\|_1\nonumber\\
 & \leq & C\,p_\ell (m(\Dt(\xi_j)\Dt)\hchi)
\label{jus}\\
 & \leq & C'\,
\sup_{\xi\in B(\be,4)}\sum_{s=0}^\ell\Dt(\xi_j)^s\,
|m^{(s)}(\Dt(\xi_j)\Dt(\xi))|\,\leq\,C''.\nonumber \Eea Thus,
raising (\ref{jus1}) to the $q^{\rm{th}}$ power  and summing we
conclude easily. \ProofEnd

\BR \Benu \item The previous proposition holds as well under the
milder hypothesis $m\in C^{n+1}(0,\infty)$ and (\ref{mih}) for
$k=0,1,\ldots,n+1$.
 This follows from the fact that the key inequality
(\ref{jus}) is actually valid with $\ell=n+1$.

\item A similar result can be proved with higher dimensional
multipliers. More precisely, given $m\in C^\infty((0,\infty)^r)$
satisfying
\[
\sup_{\xi\in(0,\infty)^r}\,
\left|\,\xi_1^{\al_1}\cdots\xi_r^{\al_r}\, \frac{\partial^{\ba }
m}{\partial\xi_1^{\al_1}\cdots\partial\xi_r^{\al_r}}
(\xi)\,\right|\,\leq \, C_{\ba},\quad\forall\;
\ba=(\al_1,\ldots,\al_r)\geq0,
\]
we can define the operator
\[
T_mf=\cF^{-1}(m(\Dt_1,\ldots,\Dt_r)\hf),\quad\mbox{if}\quad
f\in\Do.
\]
Then, an analogous proof gives $\|T_mf\|_\Bpq\leq C\|f\|_\Bpq$.

\item In the previous examples $m$ is a Fourier multiplier
belonging to $M_p$ for all $1<p<\infty$. More general examples of
multipliers for $\Bpq$ can be constructed as follows. Let
$\{m_j\}_j$ be a uniformly bounded family of multipliers in
$M_p(\SRn)$ with $\supp m_j$  contained in a fixed compact set of
$\O$. Let $m(\xi)=\sum_jm_j(g_j^{-1}\xi)$ and
$T_mf=\cF^{-1}(m\hf)$. Then $\|T_mf\|_\Bpq\leq C\|f\|_\Bpq$. In
particular, we may take $m_j=\e_j\hpsi$, where $\e_j=\pm1$, and
conclude that $m_{\bf \varepsilon}=\sum_j\e_j\psi_j$ is a
multiplier for $\Bpq$. Observe however, that letting $\e_j=1$, the
function $m_\e=\chi_\O\notin M_p(\SR^n)$ for any $p\not=2$. \Eenu
\ER

\subsection{Fourier-Laplace extensions}

It is well-known that to every distribution supported in a closed
cone $\clO$ we can associate an analytic function in the tube
domain $\Tu$ via the \emph{Fourier-Laplace integral}. More
precisely, this is given by:
\[
{\cLL} g(z)=(g,e^{i(z|\xi)})=\int_{\O}
e^{i(z|\xi)}\,g(\xi)\,d\xi,\quad z\in\Tu,
\]
which makes sense for compactly supported distributions $g$  in
$\O$, and can also be given a meaning for all $g\in\cS'(\SR^n)$
with $\supp g\subset\clO$ (see \cite[Ch. VII]{Hor}).

 In this section we wish to describe the analytic functions associated
 with (classes of) distributions in our Besov spaces $\Bpq$.
To avoid dealing with equivalence classes, it is convenient to
restrict the indices $\nu,p,q$ so that $\Bpq$ can be embedded in
the usual space of tempered distributions.

\begin{lemma}
\label{FL1} Let $\nu>0$, $1\leq p<\infty$, $1\leq q<\tqn$. Then,
for every $f\in\Sop$ with $\|f\|_\Bpq<\infty$, the series $\sum_j
f*\psi_j$ converges in the space $\cS'(\SRn)$. Moreover, the
correspondence \Beas \Bpq \;&\longrightarrow &\;\cS'(\SR^n)\\
f+\Sbop&\longmapsto &f^\sharp=\sum_j f*\psi_j\Eeas is continuous,
injective, and does not depend on the Littlewood-Paley functions
$\{\psi_j\}$.
\end{lemma}

\Proof The proof of the convergence of the series is completely
analogous to that of Lemma \ref{obser}. In fact, using the
H\"older-type inequality in (\ref{bholder}) we can write \Be
\sum_j |\lan{f*\psi_j}\phi|\leq \,C\,\|f\|_\Bpq\,
\|\phi\|_{\Bp},\quad \phi\in\cS(\SRn). \label{limBB} \Ee Now, from
$q<\tqn$  and Proposition \ref{Rsch} we obtain \Be
\|\phi\|_{\Bp}\leq
\,C\,p_\ell(\hphi)<\infty,\quad\forall\,\phi\in\cS(\SR).\label{limBB1}\Ee

Using the previous two formulas and the density of $\Do$ it is
customary to verify the last statement of the lemma. \ProofEnd

\BR {\rm From now on, whenever we restrict to the indices $\nu>0$,
$1\leq p<\infty$ and $1\leq q<\tqn$ as in the previous lemma, we
shall identify $\Bpq$ with the corresponding space $(\Bpq)^\sharp$
of tempered distributions. Observe that all these distributions
have finite order, are supported in $\clO$ and satisfy the
Littlewood-Paley decomposition $f=\sum_jf*\psi_j$.} \ER

For the purposes of this paper, we shall rather speak of
\emph{Fourier-Laplace extensions} for distributions $f$ with
spectrum in $\clO$. That is, we define
\[
{\cL}f(z)={\cLL} \hf(z)=(\hf,e^{i(z|\xi)})=\int_{\O}
e^{i(z|\xi)}\,\hf(\xi)\,d\xi,\quad z\in\Tu,
\]
which we shall call \emph{the Fourier-Laplace extension of $f$},
and which defines a holomorphic function in the tube domain $\Tu$.
For distributions $f\in\Bpq$ this takes the form \Be \cLt
f(z)=\cL(\sum_jf*\psi_j)(z)= \sum_j\cLL(\hf\hpsi_j)(z),\quad
z\in\Tu,\label{FLs} \Ee where as we shall see below, the last
series converges uniformly in compact sets of $\Tu$. We stick to
the notation $\cLt f$ to recall that we are choosing the special
representative $f^\sharp=\sum_jf*\psi_j$ from the equivalence
class $f+\Sbop$. Observe that in general $\cLt(\Sbop)=0$, while
$\cL(\Sbop)$ is not. The main result about Fourier-Laplace
extensions is the following proposition.

\begin{proposition}
\label{FL2} Let $\nu>0$, $1\leq p<\infty$, $1\leq q<\tqn$. Then
 for every distribution $f=\sum_jf*\psi_j\in\Bpq$  the series in (\ref{FLs}) converges
 uniformly on compact sets to
 a holomorphic function in $\Tu$, and moreover
\Be |\cLt f(x+iy)|\leq\,C\,\Dt(y)^{-\frac nr\frac
1p-\frac\nu{q}}\, \|f\|_\Bpq,\quad x+iy\in\Tu. \label{FL2a} \Ee In
addition, for each $y\in\O$, the distributions $\cLt f(\cdot+iy)$
satisfy
\[
\cLt(f)(\cdot+iy)=\sum_j\cLt(f)(\cdot+iy)*\psi_j,\quad\mbox{in }
\cS'(\SR^n)
\]
and \Be \|\cLt f(\cdot+iy)\|_{\Bpq}\leq C\,\|f\|_\Bpq,\Dquad
\lim_{{y\to0}\atop{y\in\O}}\|\cLt f(\cdot+iy)-f\|_\Bpq=0.
    \label{limB}
    \Ee
\end{proposition}

\Proof In the first part we shall only prove the pointwise
convergence and (\ref{FL2a}). The proof can be easily adapted to
obtain uniform convergence on compact sets. Since $\Bpq$ is
invariant by the action of $G$ as well as by translations in the
$x$ variable, we can reduce to the case $y=i\be$ and $x=0$, using
(\ref{g-action}) to prove (\ref{FL2a}) in the general case. Now,
by definition of $\cL$ and following the same steps as in
(\ref{bholder}), we can write \Beas \sum_j|\cL(f*\psi_j)(i\be)| &
= & \sum_j|\lan{f*\psi_j}{\cF^{-1}(\hchi_je^{(-\be|\cdot)})}|\\ &
\leq & C\,\|f\|_\Bpq\, \|\cF^{-1}(\chi_\O
e^{(-\be|\cdot)})\|_{\Bp}. \Eeas So,  it suffices to compute the
norm of $h=\cF^{-1}(\chi_\O
e^{-(\be|\cdot)})=\Ga_\O(\fnr)\,\Dt^{-\fnr}((\cdot+i\be)/i)$. Now,
proceeding as in (\ref{p17}) with the function $h$ we obtain the
estimate \Be |h*\psi_j(x)|\leq\,C\,\Dt^\fnr(\xi_j)\,
\frac{(1+|\xi_j|)^{2n}}
    {(1+|g_jx|^2)^n}\,e^{-\frac1\gamma(\be|\xi_j)}
\leq\,C'\,\Dt^\fnr(\xi_j)\, \frac{e^{-\frac1{2\gamma}(\be|\xi_j)}}
    {(1+|g_jx|^2)^n}.
\label{exp3} \Ee Taking $L^{p'}$-norms  and summing, we are led to
\[
\|h\|_{\Bp}\leq C\,\left[\,\int_\O
\Dt(\xi)^{\nu\frac{q'}q+\fnr\frac{q'}p}\,
e^{-c(\be|\xi)}\dxi\,\right]^\frac1{q'}.
\]
By Lemma \ref{int} this integral is finite provided
 $q<\tqn$.

\smallskip

For the second part, fix $y\in\O$ and $j$, and use the Dominated
Convergence Theorem with $\sum_k|\cL(f*\psi_k)(x+iy)|\leq C_y$
(from the first part) to write \Beas
\left(\cLt(f)(\cdot+iy)*\psi_j\right)(x) & = & \sum_k
\left(\cL(f*\psi_k)(\cdot+iy)*\psi_j\right)(x)\\
& = & \sum_k\int_{\O}\hf(\xi)\hpsi_k(\xi)e^{-(y|\xi)}\,
\hpsi_j(\xi)e^{i(x|\xi)}\,d\xi\\
& = & \int_{\O}e^{i(x+iy|\xi)}\, \hf(\xi)\hpsi_j(\xi)\,d\xi\,=\,
\cL(f*\psi_j)(x+iy). \Eeas Summing in $j$ and using again the
previous step it follows that
\[\cLt(f)(x+iy)=\sum_j\left(\cLt(f)(\cdot+iy)*\psi_j\right)(x),
\]
converging uniformly and absolutely in $x$, and hence also in
$\cS'(\SR^n)$.
\smallskip

 Let us finally prove the statements in
(\ref{limB}). First of all,
    \Bea
    \|\cLt f(\cdot+iy)\|^q_\Bpq & = &
\sum_j\Dt^{-\nu}(\xi_j)\,\|\cF^{-1}(\hf\hpsi_je^{-(y|\cdot)})\|_p^q
    \nonumber\\
     & \leq & \sum_j\Dt^{-\nu}(\xi_j)\,\|f*\psi_j\|^q_p
     \,\|\cF^{-1}(\hchi_je^{-(y|\cdot)})\|^q_1\leq
     \,C\,\|f\|^q_\Bpq,\nonumber
\Eea
 where in the last step we have estimated the $L^1$-norm by a
Schwartz seminorm \Bea \|\cF^{-1}(\hchi_je^{-(y|\cdot)})\|_1 & = &
\|\cF^{-1}(\hchi e^{-(g_jy|\cdot)})\|_1 \leq C\,p_\ell(\hchi
e^{-(g_jy|\cdot)})\nonumber\\ & & \nonumber \\
 & \leq &
C'\,(1+|g_jy|)^\ell\,e^{-\frac1\gamma(g_jy|\be)} \leq C'',
\quad\forall\,j,\,\,\,\forall\,\,y\in\O. \label{exp1} \Eea
Finally, for the convergence, use the density to write $f=g+h$,
with $g\in\cD_\O$ and $h$ with a small $\Bpq$-norm. It is well
known that $\cL g(z)$ is smooth up to the boundary with
convergence of $\cL g(\cdot+iy)$ to $g$ in the
$\cS(\SR^n)$-topology (and by  Proposition \ref{Rsch} also in the
$\Bpq$-topology). The convergence for $f$ then follows from a
standard $\e/3$ argument.

\ProofEnd

\BR\label{optimal} Let us finally remark that the index $\tqn$ in
the previous proposition is optimal. Indeed, from the duality
theorem, the continuity of the linear form $f\mapsto \cLt f(i\be)$
will imply that $\cF^{-1}(\chi_\O e^{-(\be|\cdot)})$ belongs to
$\Bp$. This continues to be the case after convolution with
$\cF^{-1}(e^{(\be|\cdot)}\hphi)$, for any  $\hphi\in
C^\infty_c(\SR^n)$ identically $1$ in a neighborhood of $0$. Thus,
it follows from  Proposition \ref{Rsch} that we must have
$q<\tqn$. \ER

\section{Bergman spaces and projectors}
\setcounter{equation}{0} \setcounter{footnote}{0}
\setcounter{figure}{0}

In this section we shall show Theorems \ref{th1} and \ref{th2}, as
well as the boundedness of the Bergman projector announced in the
introduction. Heuristically, the correspondence between
holomorphic functions $F$ in the Bergman space $\Apq$ and
distributions $f$ in $\Bpq$, is given by the Fourier-Laplace
formula:
    \[
     F(z)=\cL f(z)=\cLL\hf(z)=\int_\O e^{i\,(z|\xi)}\hf(\xi)\,d\xi, \quad
    z\in\Tu.
    \]
The distribution $f$ plays the role of a Shilov boundary value for
the holomorphic function $F$. The main result in this section is
the equivalence of norms $\|F\|_\Apq\sim\|f\|_\Bpq$, which follows
from a suitable discretization of the integral above using the
Whitney decomposition in $\S2$. Several technical estimates will
appear in this process, involving gamma integrals in $\O$ and
Littlewood-Paley inequalities as in (\ref{LP}), forcing us at some
point to assume further restrictions in the indices $\nu,p,q$. The
sharp range of parameters for the equivalence of these two norms
is still an open question, related to finer problems in Harmonic
Analysis such as restriction and cone multipliers (see the
Appendix for a further discussion on these matters).

\medskip

Before going into the proof of the theorems we need some
 preliminaries. Recall that a holomorphic function
$F\in\cH(\Tu)$ belongs to the {\it Hardy space} $H^2(\Tu)$ when

    \[
    \|F\|_{H^2}=\sup_{y\in\O}\|F(\cdot+iy)\|_{L^2(\SRn)}<\infty.
    \]
The following result is known as the \emph{Paley-Wiener Theorem}
for Hardy spaces (see Chapter III of \cite{SW}):

\begin{proposition}
\label{hardy1} A function $F\in H^2(\Tu)$ if and only if $F=\cLL
\hf$ for some $f\in L^2(\SRn)$
 with $\supp\hf\subset\clO$. In
this case, $\|F\|_{H^2}=\|f\|_{L^2(\SRn)}$.
\end{proposition}

 We will use the previous result in combination with the next one,
 whose proof is a simple modification of
 the one presented in \cite{BBPR} (see also \cite{GAR}).

\begin{proposition}
Let $\nu>\ind$ and $1\leq p,q<\infty$.
 Then, the norms of the
spaces $\Apq$ are complete. Moreover, the intersection
$H^2(\Tu)\cap A^{p,q}_\nu$ is dense in $\Apq$.
\end{proposition}

Our last preliminary  result will be the starting point in the
discretization steps to follow.

\begin{proposition}
\label{apqdiscr} Let $\nu>\ind$, $1\leq p,q<\infty$. Then, for
every lattice $\{y_j\}$ in $\O$ there exists $c>0$ such that \Be
\frac 1c\,\|F\|_\Apq\,\leq\,
\Bigl(\sum_j\Dt^\nu(y_j)\|F(\cdot+iy_j)\|_p^q\Bigr)^\frac 1q\,\leq
c\,\|F\|_\Apq,\quad \forall\, F\in\Apq(\Tu). \label{mvp2} \Ee
\end{proposition}

\Proof

 The proof relies on the
following elementary lemma:

\begin{lemma}
\label{mvp} Let  $1\leq p,q<\infty$. Then, for every
$F\in\cH(\Tu)$ and $y_0\in\O$ we have

    \Be
    \|F(\cdot+iy_0)\|_p\leq
    \,C\,\left[\int_{B(y_0,1)}\|F(\cdot+iy)\|_p^q\,\frac{dy}
    {\Dt^\fnr(y)}\right]^\frac1q,
    \label{mvp1}
    \Ee
where the constant  depends only on $p,q$.
\end{lemma}

\Proofof{Lemma \ref{mvp}}

By homogeneity, we may assume $y_0=\be$. Let us consider first the
case $p\leq q$. Then, by the mean value property for subharmonic
functions,
\[
|F(x+i\be)|^p\leq c\,\int_{B(\be,1)}\int_{|x'|\leq1}
|F(x+x'+iy)|^p\,dx'dy.
\]
Thus, integrating in $x$ and using H\"older's inequality we obtain
\[
\|F(\cdot+i\be)\|^p_p  \leq c\,\int_{B(\be,1)}
\|F(\cdot+iy)\|_p^p\,dy
  \leq
c\,\left(\int_{B(\be,1)} \|F(\cdot+iy)\|_p^q\,dy\right)^\frac{p}q.
\]
Suppose instead that $q\leq p$. Then, the mean value property
gives
\[
|F(x+i\be)|^p\leq c\,\left(\int_{B(\be,1)}\int_{|x'|\leq1}
|F(x+x'+iy)|^q\,dx'dy\right)^\frac{p}q.
\]
A new integration in $x$ and  Minkowski's inequality gives the
same result. \ProofEnd

Continuing with the proposition, we assume for simplicity that
$\{y_j\}$ is a $(\frac12,2)$-lattice, that is, $\{B_1(y_j)\}_j$
covers $\O$ with the finite intersection property. The right hand
side of (\ref{mvp2}) follows from  a discretization of the
integral on $\O$ defining $\|F\|_\Apq^q$ as in Proposition
\ref{discr}. We then use the previous lemma to conlude. For the
left hand side we just need a kind of
 converse for (\ref{mvp1}),
where as before we can assume $y_0=\be$. But this follows from the
fact that, when $F\in\Apq$,
 the function $y\to\|F(\cdot+iy)\|_p$ is monotonic in $\O$ (see,
 e.g., \cite{BBPR,GAR}), and therefore
 \[
\left[\int_{B(\be,1)}\|F(\cdot+iy)\|_p^q\,\frac{dy}
    {\Dt^\fnr(y)}\right]^\frac1q\leq C\,\|F(\cdot+ic\be)\|_p,
 \]
for some constant $c=c(\O)>0$. This  establishes the proposition.
\ProofEnd

\subsection{The proof of Theorem \ref{th1}}

We wish to show that every $F\in\Apq$ can be written as $F=\cLt f$
for some distribution $f\in\Bpq$.
 Suppose first that $F$ belongs to the dense set
 $H^2(\Tu)\cap\Apq$, so that, by Proposition \ref{hardy1},
 $F=\cLL \hf$ for some function $f\in
L^2(\SRn)$ with $\supp\hf\subset\clO$. We shall show the
inequality $\|f\|_\Bpq\leq C\,\|F\|_\Apq$. Observe that, since
$f=\sum_jf*\psi_j$ in $L^2$ (hence in $\cS'(\SR^n)$), it will
follow from this and the definition of $\cLt$ that $\cLt f=\cLL
\hf$.

\smallskip

To prove the inequality of norms, let us first denote
$y_j=\xi_j^{-1}$ the dual lattice of $\{\xi_j\}$. Then, Young's
inequality gives

    \[
    \|f*\chi_j\|_p=\|\cF^{-1}(\hf(\xi)e^{-(y_j|\xi)}
    \hchi_j(\xi)e^{(y_j|\xi)})\|_p\leq
    \|\cF^{-1}(\hf e^{-(y_j|\cdot)})\|_p\,
    \|\cF^{-1}(\hchi_je^{(y_j|\cdot)})\|_1.
    \]
Since $\xi_j^{-1}=g_j^{-1}\be$ and $g_j$ is self-adjoint, we
observe that the last factor is actually constant,
\[
\|\cF^{-1}(\hchi_je^{(y_j|\cdot)})\|_1= \|\cF^{-1}(\hchi
e^{(\be|\cdot)})\|_1=c_1<\infty.
\]
This leads to the estimate

    \Bea
    \|f\|^q_\Bpq & \leq & c\sum_j\Dt^{-\nu}(\xi_j)\,\|f*\chi_j\|_p^q\nonumber\\
     & \leq & c'\,\sum_j\Dt^{\nu}(y_j)\,\|\cF^{-1}(\hf e^{-(y_j|\cdot)})\|_p^q\nonumber\\
     & = & c'\,\sum_j\Dt^{\nu}(y_j)\,\|F(\cdot+iy_j)\|_p^q
     \,\leq\,c''\,\|F\|^q_\Apq,\label{n11}
    \Eea
where in the last inequality we have used Proposition
\ref{apqdiscr}.

For general $F\in\Apq$ one proceeds by density. Approximating
 with $\{F_m\}$ in $H^2(\Tu)\cap\Apq$,
 we obtain a corresponding
sequence of functions $f_m\in L^2(\SRn)$, which by (\ref{n11}) is
a Cauchy sequence in $\Bpq$. Then, by completeness of this space
and Lemma \ref{FL1}, $f_m$ converges (in $\Bpq$ and in $\cS'$) to
a distribution $f\in \Bpq$ such that $f=\sum_jf*\psi_j$. Moreover,
$\|f\|_\Bpq\leq C\,\|F\|_\Apq$. It remains to prove that
$F=\cLt(f)$, for which we can use the continuity of the pointwise
evaluation functional $f\mapsto\cLt f(z)$ in (\ref{FL2a}). Indeed,
for each $z\in\Tu$ we have
\[
\cLt f(z)=\lim_{m\to\infty}\cLt f_m(z)=\lim_{m\to\infty}\cLL
\hf_m(z)=\lim_{m\to\infty}F_m(z)=F(z).
\]
Finally, the convergence \[
\lim_{{y\to0}\atop{y\in\O}}F(\cdot+iy)=f, \quad\mbox{in }\quad
\Bpq\mand\cS'(\SR^n),\] is just a consequence of Proposition
\ref{FL2} and Lemma \ref{FL1}. This completes the proof of the
theorem. \ProofEnd

\subsection{The proof of Theorem \ref{th2}}

We start with a preliminary result which is the dual version of
(\ref{LP}),  quoted in the introduction.

\begin{lemma}
\label{Lplemma} Let $1\leq p<\infty$ and let $1\leq s\leq
p_\sharp=\min\{p,p'\}$. Then there exists a constant $C$ such
that, for every sequence of functions  $f_j\in L^p(\SRn)$
satisfying $\supp\hf_j\subset B_2(\xi_j)$,
we have the inequality
 \Be \|\sum_jf_j\|_p
\leq\,C\Bigl(\sum_j\|f_j\|_p^s\Bigr)^\frac1s . \label{LP1} \Ee
\end{lemma}

\Proof  It is sufficient to prove the stronger
 inequality $$\|\sum_jf_j*\chi_j\|_p
\leq\,\Bigl(\sum_j\|f_j\|_p^s\Bigr)^\frac1s,$$ valid for all
sequences of functions in $L^p$. The  $\chi_j$ are chosen as in
$\S3.1$ with their $L^1$ norm uniformly bounded, and their Fourier
transform supported in  $B_4(\xi_j)$ and identically $1$ on the
ball $B_2(\xi_j)$. For this last inequality, the proof is
immediate when $s=1$ by Minkowski's inequality, as well as for
$s=p=2$ by the finite intersection property of the balls. We
interpolate between these two cases to conclude. \ProofEnd

\BR \label{Rmain}{
 Our proof for Theorem \ref{th2} depends
directly on the previous simple inequality, in which unfortunately
the best exponent $s$ for each fixed $p$ seems not to be known
(ideally, $s=2$ would be the best possible). The reader can track
down in our proof below that any improvement over $s=p_\sharp$,
for a fixed $p$, will end up in the validity of the theorem for
all $q<sq_\nu$ and the same $p$. For this reason we state below a
result, where a more general inequality than (\ref{LP1}) is
assumed to hold true. We shall discuss in the Appendix the
validity of such inequalities in the case of light-cones.}\ER

\medskip

\begin{proposition} \label{thR1}

\noindent Let $\nu>\ind$, $1\leq p,s<\infty$. Assume there exist
numbers $\mu,\dt\geq 0$ and a constant $C=C(\mu,\dt)>0$ such that
  \Be \|\sum_jf _j\|_p\leq C \,\Bigl[\,\sum_j\Dt^{-\mu}(\xi_j)e^{\delta(\xi_j|\be)}
\|f_j\|_p^s\,\Bigr]^\frac1s \label{R5} \Ee   holds for every
finite sequence $\{f_j\}\subset L^p(\SRn)$ satisfying
$\supp\hf_j\subset B_2(\xi_j)$. Then, for every index
\Be
q<\min\left\{\,s\tfrac{q_\nu}{q_\mu}\,,
\,s\,\tfrac{\nu-(\fnr-1)}\mu\,,\,\tqn\,\right\}, \label{R6} \Ee
and for every distribution $f$ with $\|f\|_\Bpq<\infty$, the
function $F=\cL (\sum_jf*\psi_j)$ belongs to $\Apq$, and moreover,
\[
\|F\|_\Apq\,\lesssim\,
 \|f\|_{B^{p,q}_\nu}.
\]
\end{proposition}

Using Lemma \ref{Lplemma}, it is clear that Proposition \ref{thR1}
suffices to establish Theorem \ref{th2}. Indeed, the lemma implies
the assumption of the proposition with $\mu=\dt=0$ and
$s=p_\sharp=\min\{p,p'\}$. In this case $q_\mu=1$, so if we assume
the validity of the proposition the condition on $q$ simplifies
into $q<\min\{sq_\nu,\tqn\}=p_\sharp q_\nu$, as stated in Theorem
\ref{th2}.

\medskip

We turn now to the proof of Proposition \ref{thR1}. We shall show
that, for every $f\in\Do$, the function $F(z):=\cL
f(z)=\cLL\hf(z)$ belongs to $\Apq(\Tu)$, and $\|F\|_\Apq\leq
C\,\|f\|_\Bpq$. This will be enough to conclude, since in the
general case one can proceed by density. We will need an
intermediate result, which we will comment later.

\begin{lemma}
\label{fLR} Let $1\leq p,s<\infty$, and assume that (\ref{R5})
holds for some $\mu, \dt\geq 0$. Then, for every $f\in\Do$ and
$y\in\O$, the function $F(\cdot +iy)=\cL f(\cdot+iy)$ belongs to
$L^p(\SR^n)$. Moreover, \Be \|F(\cdot +i y)\|_p\,\lesssim\,
 \Dt^{-\frac{\mu}s}(y)\|f\|_{B^{p,s}_\mu}
\label{R3}  \Ee with constants independent of $f$ or $y\in\O$.
\end{lemma}

\Proof By homogeneity (see Lemma \ref{besov-norm-lem}), it is
sufficient to prove (\ref{R3}) when $y=\eta\be$, for some fixed
$\eta>0$ to be chosen below. Now, given $f\in\Do$, let us define
$g\in\Do$ by $\hg=\hf e^{-\eta(\be|\cdot)}$. Then, applying the
assumption to $g=\sum_j g*\psi_j$, we obtain \Beas \|g\|_p & = &
\|F(\cdot+i\eta\be)\|_p\;\; \lesssim\;
\Bigl(\,\sum_j\Dt^{-\mu}(\xi_j)e^{\delta(\xi_j |
\be)}\|\cF^{-1}(\hf\hat\psi_j
e^{-\eta(\be|\cdot)})\|_p^s,\;\Bigr)^\frac1s\\
  & \lesssim &
\Bigl(\,\sum_j\Dt^{-\mu}(\xi_j)e^{\delta(\xi_j |
\be)}\|f*\psi_j\|_p^s\,
\|\cF^{-1}( e^{-\eta(\be|\cdot)}\hchi_j)\|_1^s\,\Bigr)^\frac1s\\
  & \lesssim &
\Bigl(\,\sum_j\Dt^{-\mu}(\xi_j)\|f*\psi_j\|_p^s\,\Bigr)^\frac1s.
\Eeas In the last step we have used the fact, from (\ref{exp1}),
that $\|\cF^{-1}(e^{-\eta(\be|\cdot)}\hchi_j)\|_1$ can be bounded,
up to some constant, by $e^{-\gamma\eta (\xi_j|\be)}$. We conclude
by choosing any $\eta$ larger than $\delta/\gamma$. \ProofEnd

\medskip

Let us go back to the proof of the proposition. Given $f\in\Do$,
and $F(z):=\cL f(z)=\cLL\hf(z)$, the previous lemma applied to
$\cF^{-1}(\hf e^{-(y|\be)})$ gives us \Beas \|F(\cdot+i2y)\|_p
  & \lesssim &\Dt^{-\frac \mu{s}}(y)\,
\Bigl[\,\sum_j\Dt^{-\mu}(\xi_j)\,
\|\cF^{-1}(\hf\hpsi_j e^{-(y|\cdot)})\|_p^s\,\Bigr]^\frac1s\\
& \lesssim &\Dt^{-\frac \mu{s}}(y)\,
\Bigl[\,\sum_j\Dt^{-\mu}(\xi_j)\,
e^{-\gamma(y|\xi_j)}\,\|f*\psi_j\|_p^s\,\Bigr]^\frac1s, \Eeas
where once again we have used (\ref{exp1}). Thus, \Beas
I & := & \int_\O\|F(\cdot+iy)\|_p^q\,\Dt^{\nu-\fnr}(y)dy\\
  & \lesssim &
\int_\O \Dt^{-\mu \frac qs}(y)\, \left(\sum_j\Dt^{-\mu}(\xi_j)
e^{-\gamma(y|\xi_j)}\|f*\psi_j\|_p^s\right)^\frac
qs\,\Dt^{\nu-\fnr}(y)dy. \Eeas  When $q\leq s$ we directly
conclude
\[
I\,\lesssim\,\sum_j \Dt^{-\mu \frac qs}(\xi_j)\|f*\psi_j\|_p^q\,
\int_\O \Dt^{-\mu\frac qs}(y)\,
e^{-\gamma'(y|\xi_j)}\,\Dt^{\nu-\fnr}(y)dy,
\]
where the gamma integral equals a multiple of
$\Dt(\xi_j)^{\mu\frac qs-\nu}$, whenever $\nu-\mu\frac qs
>\fnr-1$. This leads to one of the conditions stated in
(\ref{R6}).

Suppose now that $q>s$. Then we multiply and divide the summands
by $\Dt_{\bt}(\xi_j)$, for some multi-index
$\bt=(t_1,\ldots,t_r)\in\SR^r$ to be chosen below. After applying
H\"older's inequality, we obtain \Beas I & \lesssim &
\int_\O\Dt^{-\mu \frac qs}(y)\, \Bigl[\,\sum_j\Dt^{-\mu \frac
qs}(\xi_j)e^{-\gamma(y|\xi_j)} \|f*\psi_j\|_p^q\,
  \Dt_{-\bt\frac qs}(\xi_j)\,\Bigr]\, \\
& & \Pquad \times\;
   \Bigl[\,\sum_j\Dt_{\bt(\frac qs)'}(\xi_j)
  e^{-\gamma(y|\xi_j)}\,\Bigr]
^{\frac{q/s}{(q/s)'}}\dmes. \Eeas According to Proposition
\ref{discr}, the last bracket can be transformed into a gamma
integral, which in order to be finite requires the following
condition in the indices
\[
t_j(q/s)'> \,(j-1)\frac{n/r-1}{r-1},\quad\mbox{for all}\,\,\,
  j=1,\ldots,r.
\]
Thus, replacing the expression inside the brackets by a multiple
of $\Dt^*_{-\bt^*(q/s)'}(y)$, we have
 \Beas I  &
\lesssim & \sum_j\Dt^{-\mu \frac qs}(\xi_j) \,\|f*\psi_j\|_p^q\,
  \Dt_{-\bt\frac qs}(\xi_j)\,
\int_\O \Dt^*_{-(\bt^*+\mu)\frac qs}(y)e^{-\gamma(y|\xi_j)}\dmes\\
   & \lesssim &
\sum_j\Dt^{-\nu}(\xi_j)\, \|f*\psi_j\|_p^q. \Eeas In the last step
we have computed again the gamma integral, which gives a finite
multiple of $\Dt_{(\bt+\mu)\frac qs-\nu}(\xi_j)$ if we impose the
condition in the indices:
\[
-\frac qs(t_{r-(j-1)}+\mu)+\nu\,>
\,(j-1)\frac{n/r-1}{r-1},\quad\mbox{for all}\,\,\,
  j=1,\ldots,r.
\]
Therefore, the $t_j$'s must be chosen so that:
     \[
     \frac1{(q/s)'}\,\frac{j-1}{r-1}\,
\left({\Ts\frac{n}r}-1\right) \,<\,t_j\,<\, \frac1{q/s}\,
\left(\nu-\frac{r-j}{r-1}\, ({\Ts\frac{n}r}-1)\,\right)-\mu, \quad
j=1,\ldots,r.
     \]
Using $\frac1{(q/s)'}=1-\frac1{q/s}$, we see that this is only
possible when
\[
\frac{j-1}{r-1}\,\bigl({\Ts\frac{n}r}-1\bigr)+\mu \,<\,
\frac1{q/s}\,\biggl(\nu-(\ind)+2\,\frac{j-1}{r-1}\,(\ind)\,\biggr),\quad
j=1,\ldots,r.
\]
Solving for $q/s$, this forces us to  have
\[
\frac{q}s<\, \min_{1\leq j\leq r} {{\nu-(\ind)+
2\,\frac{j-1}{r-1}(\frac{n}r-1)}
\over{\mu+\frac{j-1}{r-1}(\frac{n}r-1)}}= \left\{\Ba{lll}
q_\nu/q_\mu, & & \nu> 2\mu+\ind;\\
& & \\
\frac{\nu-(\ind)}\mu, & & \nu\leq 2\mu+\ind.\Ea\right.
\]
This is precisely the range of $\nu$ and $q$ assumed in (\ref{R6})
 so the theorem is completely proved.
 \ProofEnd

We conclude this subsection with some equivalent versions of the
sufficient condition (\ref{R5}) when $\mu>0$.  These are
convenient expressions which we shall relate in the Appendix with
known inequalities in light-cones.

\begin{proposition}\label{equivalent}
Let $1\leq p,s <\infty$ and $\mu>0$. Then the following properties
are equivalent: \begin {itemize} \item [1.] There exists
$\delta>0$, and a constant $C_\delta>0$ such that  \Be
\bigl\|\sum_jf_j\bigr\|_p\, \leq \,C_\delta\,
\Bigl[\,\sum_j\Dt^{-\mu}(\xi_j)\,e^{\delta(\xi_j|\be)}\,\|f_j\|_p^s\,\Bigr]^\frac1s,
 \label{R14}
 \Ee
for every finite sequence $\{f_j\}$ in $L^p(\SRn)$ satisfying
$\supp\hf_j\subset B_2(\xi_j)$.

\item [2.] There exists a constant $C>0$ such that \Be
\|\sum_jf_j\|_p\, \leq \,C\,
\Bigl[\,\sum_j\Dt^{-\mu}(\xi_j)\|f_j\|_p^s\,\Bigr]^\frac1s,
 \label{R13}
 \Ee
 for every finite sequence $\{f_j\}$ in $L^p(\SRn)$ satisfying
$\supp\hf_j\subset B_2(\xi_j)\cap H_0$, where $H_0$ is the band in
between two hyperplanes $H_0=\{\frac12<(\be|\xi)<2\}$.

 \item [3.] There
exists a constant $C>0$ such that
 \Be
 \|\cL f(\cdot+iy)\|_p
\leq\, C \, \Dt^{-\frac{\mu}s}(y)\,\|f\|_{B^{p,s}_\mu},\quad
\forall\;f\in\Do,\quad\forall \,y\in\O. \label{R12}  \Ee
\end{itemize}
\end{proposition}

\Proof We have already proved in Lemma \ref{fLR} that
(\ref{R14})$\Rightarrow$(\ref{R12}). To prove that
(\ref{R12})$\Rightarrow$(\ref{R13}), take a corresponding sequence
$\{f_j\}$ as in (\ref{R13}). Define the functions
$\hg_j=e^{(\be|\cdot)}\hf_j$, and apply the inequality (\ref{R12})
for $y=\be$ to the function $g=\sum_jg_j\in\Do$. Then,
$$\|\sum_j f_j\|_p= \|\cF^{-1}(\hg e^{-(\be|\cdot)})\|_p \leq
\,C \,\Bigl[\,\sum_k\Dt^{-\mu}(\xi_k)\,\|\cF^{-1}(\hg\hat\psi_k)
\|_p^s\,\Bigr]^\frac1s.$$ Using the finite intersection property,
and the fact that the $L^p$-norms of $f_j$ and $g_j$ are
comparable, we obtain easily the right hand side of (\ref{R13}).

It remains to prove that (\ref{R13})$\Rightarrow$(\ref{R14}). To
do this, we are going to slice the cone with hyperplanes, and then
apply an scaled version of (\ref{R13}) to the restrictions of
$\sum_jf_j$ to the bands
\[
H_k=\{2^{k-1}<(\xi|\be)<2^{k+1}\},\quad k\in\SZ.
\]
To do this argument precise, we select a sequence of smooth
1-variable functions $\{\rho_k\}$ so that $\supp\rho_k\subset
(2^{k-1},2^{k+1})$ and $\sum_{k\in\SZ}\rho_k\equiv1$ in
$(0,\infty)$. We let $\hf_{j,k}(\xi)=\rho_k((\xi|\be))\hf_j(\xi)$,
so that \Be \supp\hf_{j,k}\subset B_2(\xi_j)\cap H_k\mand
\|f_{j,k}\|_p\leq C\,\|f_j\|_p. \label{aux3}\Ee By Minkowski's
inequality we can write
$$I=\,\bigl\|\sum_jf_j\bigr\|_p\, \leq
 \,\sum_{k \in\SZ}\,\|\sum_{j\in J_k}f_{j,k}\|_p\,=:\,\sum_{k\in\SZ} \|F_k\|_p, $$
where the sets of indices  $J_k$ are defined so that
$B_2(\xi_j)\cap H_k\not=\emptyset$. In order to estimate the norm
of $F_k$, we must first perform a dilation by $\dt=2^{-k}$ so
that, replacing $F_{k}$ with $F^{(\dt)}_{k}=\dt^\frac
np\,F_{k}(\dt\,\cdot)$, we do not change the $L^p$-norms and the
Fourier transform is now supported in $H_0$.  Thus, we are in
conditions of applying (\ref{R13}):
\[
\|F_k\|_p=\|\sum_\ell F^{(\dt)}_k*\psi_\ell\|_p\leq\,C\,
\Bigl[\,\sum_\ell\Dt^{-\mu}(\xi_\ell)\,\|F^{(\dt)}_k*\psi_\ell\|_p^s\,\Bigr]^\frac1s.
\]
Now, observe that each $f^{(\dt)}_{j,k}$ has Fourier transform
supported in $B_2(\delta\xi_j)$, and $\{\dt\xi_j\}$ is still a
$(\frac12,2)$-lattice in the cone. Thus, by the finite
intersection property, the set of indices $j$ for which
$B_2(\delta\xi_j)$ intersects a fixed set $B_2(\xi_\ell)$ has at
most $N=N(\O)$ elements, independently of $\ell$ and $\dt$. Thus,
\[
\|F_k\|_p^s\lesssim
\sum_\ell\Dt^{-\mu}(\xi_\ell)\sum_j\|f^{(\dt)}_{j,k}*\psi_\ell\|_p^s.
\]
Now, changing the order of sums, restricting $\ell$ to the bounded
set of indices
\[
\tJ_j=\{\ell\mid B_2(\xi_\ell)\cap
B_2(\delta\xi_j)\not=\emptyset\},
\]
and using that $\Dt(\xi_\ell)\sim\Dt(\dt\xi_j)$ for such indices,
we obtain \Beas \|F_k\|_p^s & \lesssim &
\sum_j\Dt^{-\mu}(\dt\xi_j)\sum_{\ell\in\tJ_j}\|f^{(\dt)}_{j,k}\|_p^s\|\psi_\ell\|_1^s\\
& \lesssim & \dt^{-\mu
r}\,\sum_j\Dt^{-\mu}(\xi_j)\|f_{j,k}\|_p^s\,  \lesssim \, 2^{k\mu
r}\,\sum_{j\in J_k}\Dt^{-\mu}(\xi_j)\|f_j\|_p^s, \Eeas where in
the last step we have also used (\ref{aux3}). Thus, raising to the
power $1/s$ and summing in $k$, we have shown that
\[
I\lesssim\sum_{k\in\SZ}\,\Bigl[\,\sum_{j\in
J_k}\Dt^{-\mu}(\xi_j)\|f_j\|_p^s\,\Bigr]^\frac1s\,2^{\frac{k\mu
r}s}.
\]
Multiplying and dividing by $e^{2^k}$, and applying H\"older we
obtain
\[
I\lesssim \Bigl[\,\sum_{k\in\SZ}\,\sum_{j\in
J_k}\Dt^{-\mu}(\xi_j)\,e^{s2^k}\,\|f_j\|_p^s\,\Bigr]^\frac1s\,
\Bigl[\,\sum_{k\in\SZ}2^{k\mu
rs'/s}\,e^{-s'2^k}\,\Bigr]^\frac1{s'}.
\]
The last term is a finite constant when $\mu>0$, while in the
first factor we can replace $e^{s2^k}$ by $e^{\eta(\xi_j|\be)}$,
for a sufficiently large $\eta$. To conclude it remains only to
show that, for each fixed $j$, the set of all $k\in\SZ$ such that
$H_k$ intersects $B_2(\xi_j)$ contains at most $N=N(\O)$ elements.
To see this observe, from Lemma \ref{l6}, that each such $k$ must
satisfy
\[
\tfrac{2^{k-1}}\gamma<(\xi_j|\be)<\gamma\,2^{k+1},
\]
or equivalently
$\frac1{2\gamma}\,(\xi_j|\be)<2^k<2\gamma\,(\xi_j|\be)$. Taking
logarithms we see that this is only possible for a constant number
of such $k$'s. The proof of Proposition \ref{equivalent} is then
complete. \ProofEnd

\subsection{Boundedness of Bergman projectors in $\Lpq$}

In this section we shall prove that the norm equivalence between
$\Apq$ and $\Bpq$ is equivalent to the boundedness of the Bergman
projector $P_\nu$ in  $\Lpq$ spaces.

\smallskip

Recall that the Bergman projector $P_\nu$ is defined for functions
$F\in\Ld(\Tu)$ as \Be P_\nu
F(x+iy)=\int_\SRn\int_{\O}B_\nu(x-u+i(y+v)) F(u+iv)
\,\Dt(v)^{\nu-\fnr}\,dvdu, \label{proj} \Ee where the  Bergman
kernel has the well-known expression \Be B_\nu(z-{\Ol{w}})=
d(\nu)\,\Dt^{-(\nu+\fnr)}((z-{\Ol{w}})/i)=c_\nu\int_\O
e^{i(z-{\Ol{w}}|\xi)}\,\Dt(\xi)^\nu\,d\xi,\quad z,w\in\Tu,
\label{berg} \Ee for some positive constants $c_\nu,d(\nu)$ (see,
e.g., Chapter XIII of \cite{FK}).
%
It is clear that $P_\nu F(z)$ defines a holomorphic function in
$\Tu$ whenever the integral in (\ref{proj}) converges absolutely.
The following lemma shows that this is the case exactly when
$F\in\Lpq$ and $q<\tqn$. This elementary fact also gives us a
trivial range of unboundedness for $P_\nu$ (see \cite{BT}).

\begin{lemma}
\label{basic} Let $\nu>\ind$ and $1\leq p<\infty$. Then \Be
B_\nu(z+i\be)\in L^{p',q'}_\nu(\Tu)\quad\Longleftrightarrow
 \quad
q<\tqn:=\frac{\nu+\fnr-1}{(\fnr\frac1{p'}-1)_+}, \label{Bnu} \Ee
Moreover, if $q\leq \tqn'$ or $q\geq\tqn$ then $P_\nu$ does not
admit bounded extensions into $\Lpq$.
\end{lemma}
\Proof The first statement is an elementary application of Lemma
\ref{ial} to the formula in (\ref{berg}). For the second statement
test with $F(z)=\Dt^{-\nu+\fnr}(\Im z)\chi_{Q(i\be)}(z)$, where
$Q(i\be)$ is a closed polydisk in $\Tu$ centered at $i\be$. Then
\[
P_\nu F(z)=c_n\,B_\nu(z+i\be), \quad z\in\Tu,
\]
by the mean value property for (anti)-holomorphic functions.
Therefore, if $q\geq \tqn$, $P_\nu$ cannot be bounded into
$L^{p',q'}_\nu$, and by self-adjointness neither into $\Lpq$.
\ProofEnd

We pass now to the study of boundedness of $P_\nu$ in $\Lpq$ when
$\tqn'<q<\tqn$. This is a difficult open question for which only
partial results are known (see \cite{BBPR} for the light-cone, and
\cite{BBG} for the simpler case $L^{2,q}_\nu$). We prove here the
following equivalence between this problem and the kind of
estimates for Fourier-Laplace integrals that we have considered.
We will see that it is an easy  consequence of our previous study.
\begin{theorem}
\label{necsuff} Let $\nu>\ind$ and $2\leq q<\tqn$. Then, the
Bergman projector $P_\nu$ is bounded in $\Lpq$ if and only if
there exists a constant $C$ such that \Be \|\cLL \hf\|_{\Lpq}\leq
C \|f\|_{\Bpq},\quad f\in\Do. \label{necsuff2} \Ee
\end{theorem}

\Proofof{\bf   the necessity condition} Let $2\leq q<\tqn$ and
assume that the projector $P_\nu$ is bounded in $\Lpq$. We want to
compute $\|\cLL \hf\|_{\Lpq}$ for $f\in\Do$. It is sufficient to
test it on functions in $L^{p',q'}_\nu \cap L^2_\nu$. Moreover,
since $\cLL\hf\in\Apq$ and the projection is self-adjoint (hence,
bounded in $L^{p',q'}_\nu $), we can as well test it on functions
which are in $A^{p',q'}_\nu \cap A^2_\nu$.  Such functions may be
written as $\cLt g=\cLL\hg$, with $g=\sum_jg*\psi_j\in
B^{p',q'}_\nu$, and we know from Theorem \ref{th2} that, for this
range of exponents, the norm of $\cLL \hg$ in $A^{p',q'}_\nu$ is
equivalent to the norm of $g$ in $B^{p',q'}_\nu$. So it is
sufficient to prove that
$$\left|\int_{\SRn}\int_\O \cLL\hf (x+iy)\overline{\cLL\hg
(x+iy)}\Dt(y)^{\nu-\frac nr} dy\,dx\right| \leq C
\|f\|_{\Bpq}\|g\|_{B^{p',q'}_\nu},$$ for some constant $C$ which
does not depend on $f$ and $g$. Using the Paley-Wiener Theorem for
$\Ad$ (see, e.g., \cite[p. 260]{FK}), we know that the left hand
side is equal to
$$\left|\,\int_\O \widehat f (\xi) \overline{\widehat g (\xi)}\frac
{d\xi}{\Dt(\xi)^{\nu}}\,\right|.$$ Then, Plancherel's Theorem and
the definition of the fractional powers of $\Box$ tell us that
this is the usual duality pairing between $f$ and $\Box^{-\nu}g$.
As a consequence, it is bounded by
$$C\|f\|_{\Bpq}\|\Box^{-\nu}g\|_{B^{p',q'}_{-\nu q'/q}}.$$
To conclude, we  use Proposition \ref{ib}, which gives the
equivalence of the norm  $\|\Box^{-\nu}g\|_{B^{p',q'}_{-\nu
q'/q}}$ with the norm $\|g\|_{B^{p',q'}_\nu}$. This finishes the
proof of this direction. \ProofEnd

\smallskip

For the other direction, we will prove a little more. We will show
 that $P_\nu$ is always bounded from $\Lpq$ into a new holomorphic function space
$\cBpq:=\cLt(\Bpq)$, consisting of Fourier-Laplace transforms of
distributions in $\Bpq$. Then one concludes easily from there and
the identification $\cBpq=\Apq$ in the previous section. To make
these statements precise we begin with the definition of $\cBpq$.

\begin{definition}
Given $\nu>0$, $1\leq p<\infty$ and $1\leq q<\tqn$, we define the
holomorphic function space
\[
\cBpq(\Tu):=\left\{F=\cLt f=\sum_j \cLL(\hf\hpsi_j)\mid
f\in\Sop\;\mbox{with}\;\|f\|_\Bpq<\infty\right\},
\]
endowed with the norm $\|F\|_\cBpq=\|f\|_\Bpq$.
\end{definition}

By Proposition \ref{FL2}, it follows that $\cLt\colon\Bpq\to\cBpq$
is an isomorphism of Banach spaces, $\cBpq$ is continuously
embedded in $\cH(\Tu)$ and its functions satisfy the inequality
\[
|F(x+iy)|\leq\,C\,\Dt(y)^{-\frac nr\frac 1p-\frac\nu{q}}\,
\|f\|_\Bpq,\quad x+iy\in\Tu.\] In this context, Theorems \ref{th1}
and \ref{th2} can be written as well as
\[
\Apq\subset\cBpq\mbox{ when }1\leq q<\tqn,\mand \Apq=\cBpq\mbox{
when }1\leq q<q_{\nu,p}. \] Observe also that $\Apq$ is a dense
subspace of $\cBpq$ (since $\cL(\Do)\subset\Apq$), but in general
is not closed. In fact, examples in the next subsection show that
the inclusion is strict whenever $q\geq\min\{2,p\}q_\nu$. We now
prove the announced statement, which allows to conclude for the
proof of Theorem \ref{necsuff}.

\begin{proposition}
\label{main2} Let $\nu>\ind$, $1\leq p<\infty$ and
 $2\leq q<\tqn$. Then, with the previous notation, $P_\nu$ extends as a bounded
operator from $\Lpq$ into $\cBpq$. That is, for every $F\in\Lpq$
there exists $g\in\Bpq$ such that $P_\nu F=\cLt g$ with \Be
\|P_\nu F\|_{\cBpq}=\|g\|_\Bpq\leq C\,\|F\|_\Lpq,\quad
F\in\Lpq(\Tu). \label{thp1}\Ee
\end{proposition}

\BR Observe that for $F\in\Ad\cap\Apq$, $P_\nu F=F=\cLt f$, and
therefore, (\ref{thp1}) actually generalizes the inequality
$\|f\|_\Bpq\leq C\,\|F\|_\Apq$ in Theorem \ref{th1}. \ER

\Proof It is enough to show (\ref{thp1}) for functions $F$ in the
dense set $\Ld\cap\Lpq$, proceeding otherwise as in the last part
of $\S4.1$. Since $P_\nu$ is a projector, for such functions we
will have $P_\nu F\in\Ad$, and therefore, by the Paley-Wiener
Theorem for $\Ad$  there exists a unique function $\hg\in
L^2(\O;\Dt^{-\nu}(\xi)d\xi)$ such that \Be P_\nu
F(z)=\cLL\hg(z)=\int_\O e^{i(x+iy|\xi)}\hg(\xi)\,d\xi,
 \quad z=x+iy\in\Tu.
 \label{FLt2}
 \Ee
Observe that $g=\sum_jg*\psi_j$ in $\cS'(\SR^n)$ (since
$\hg=\sum_j\hg\hpsi_j$ in $L^2(\O;\Dt^{-\nu}(\xi)d\xi)$), and
therefore $P_\nu F=\cLL\hg=\cLt g$. Thus, to  prove $P_\nu
F\in\cBpq$ we just need to bound $\|g\|_\Bpq$. From the duality of
Besov spaces (Proposition \ref{bpqdual}), it follows that it is
sufficient to prove that
$$|\lan g\phi |\leq C \|F\|_{\Lpq}\|\phi\|_{B^{p',q'}_{-\nu q'/q}},\quad\phi\in\Do.$$
As in the previous proof, we use the fact that
\begin{eqnarray*}
   \lan g\phi  & =\int_{\SRn}\int_\O \overline{\cLL\hg (x+iy)}\,{\cLL\hh
(x+iy)}\Dt(y)^{\nu-\frac nr} dy\,dx \\
     & \int_{\SRn}\int_\O \overline{F(x+iy)}\,{\cLL\hh
(x+iy)}\Dt(y)^{\nu-\frac nr} dy\,dx,
  \end{eqnarray*}
with $h=\Box^{\nu}\phi $. So,
$$|\lan g\phi|\leq C \|F\|_{\Lpq}\|\cL h\|_{A^{p',q'}_{\nu}}.$$
To conclude, we use Theorem \ref{th2} applied to $h$, and
Proposition
 \ref{ib} as before to have the equivalence of the norm
$\|\phi\|_{B^{p',q'}_{-\nu q'/q}}$ with the norm
$\|h\|_{B^{p',q'}_\nu}$. \ProofEnd

As a corollary, we can use the previous section to extend the
range of exponents for which the Bergman projector is bounded.

\begin{corollary}
If $\nu>\ind$, $1\leq p<\infty$ and $q'_{\nu,p}<q<q_{\nu,p}$, then
$P_\nu$ admits a bounded extension to $\Lpq$. That is, there
exists a constant $C>0$ such that
\[
\|P_\nu F\|_\Lpq\leq\,C\,\|F\|_\Lpq, \quad\forall\,\,F\in\Lpq.
\]
\end{corollary}

\subsection{Necessary conditions}

In this section we shall construct counter-examples to Theorem
\ref{th2}, for some values of $\nu, p,q$ above the critical
indices. That is, we shall show that the space $\cBpq$, of
Fourier-Laplace extensions of $\Bpq$, cannot be embedded into
$\Apq$. Our examples are actually stronger and show that $\cBpq$
cannot be even embedded into $A^{p,\infty}_\nu$, that is, there
\emph{does not} exist a constant $C$ such that, for all $f\in
\Bpq$, one has the inequality \Be \int_{\SRn}|F(x+i\be)|^p\,dx\leq
C\,\|f\|_\Bpq^p,\quad \mbox{for}\quad F= \cLt(f).
 \label{cannot}\Ee
 Observe that this indeed contradicts Theorem \ref{th2} since, by Lemma \ref{mvp},
 the integral on the left hand side is always smaller than $\|F\|_\Apq$.
 Our results are stated in the following proposition.
 \begin{proposition}\label{nec} Let $\nu>\ind$, $1\leq p<\infty$ and
$1\leq q<\tqn$. Then, there cannot exist a constant $C$ such that
the inequality (\ref{cannot}) is valid for all $f\in \Bpq$ in the
two following cases: \Benu \item[{\rm(a)}] $1\leq p\leq 2$ and
$q\geq q_{\nu,p}$;

\item[{\rm(b)}] $2<p<\infty$ and $q\geq \min\{2 q_\nu,\tqn\}$.
\Eenu
\end{proposition}
\Proof
 We shall use a different method for (a) and (b).
 The first one is based on an explicit holomorphic function, and the second on a
 Rademacher argument with Littlewood-Paley inequalities.
 For the first part,
 we shall find  $F\in\cH(\Tu)$ such that
\[
\int_{\SRn}|F(x+i\be)|^p\,dx=\infty\mand \|\Box
F\|_{L^{p,q}_{\nu+q}}<\infty,
\]
for all $q\geq pq_\nu$ (with $q<\tqn$). This gives a contradiction
with (\ref{cannot}), since in case that held, we would conclude
\[
\|F(\cdot+i\be)\|_p\lesssim\,\|f\|_\Bpq\sim\|\Box
f\|_{B^{p,q}_{\nu+q}}\leq C\,\|\cLt(\Box
f)\|_{L^{p,q}_{\nu+q}}=C\|\Box F\|_{L^{p,q}_{\nu+q}}<\infty.
\]
At this point, an example involving the $\Box$ operator may seem
a bit cumbersome, but it is actually quite natural since the
boundedness of the Bergman projection turns out to be equivalent
to the existence of \emph{generalized Hardy inequalities}. For
more on this direction we refer to \cite{BBPR} (in the case of the
light-cone), and to the survey paper \cite{B}. Our specific
example will be the holomorphic function
\[
F(z)=\Dt((z+i\be)/i)^{-\al}\,(1+\log\Dt((z+i\be)/i))^{-\frac1p},\quad
z\in\Tu,
\]
where we choose $\al=(\frac{2n}r-1)/p$. We are using the standard
convention
$\log\Dt(z/i)=\sum_{j=1}^r\log[\frac{\Dt_j}{\Dt_{j-1}}(z/i)]$,
since in this case $\re[\frac{\Dt_j}{\Dt_{j-1}}(z/i)]>0$, for each
$z\in\Tu$ and $j=1,\ldots,r$ (see, e.g., the discussion in
\cite[$\S7$]{GAR}). We also remark that, since
$|\Dt((z+i\be)/i)|\geq \Dt(y+\be)> 1$, the expression under the
power $\frac1p$ has positive real part, and defines a holomorphic
function.

To compute the first integral we estimate the denominator of
$F(z)$ using the elementary facts \[|1+\log\Dt((x+i\be)/i)|\leq
r{\Ts\frac\pi2}+1+\log|\Dt(x+i\be)|\mand |\Dt(x+i\be)|\leq
\Dt(x+\be),
\]
when $x\in\O$.  This leads to the expression
\[
\int_{\SRn}|F(x+i\be)|^p\,dx\geq C\,\int_\O \frac{dx}
{|\Dt(x+\be)|^{\frac{2n}r-1}\,(1+\log\Dt(x+\be))}=\infty,
\]
 by Lemma \ref{log2}.

For the second integral we first calculate $\Box F(z)$. Observe
that around each $z_0\in\Tu$ there is a neighborhood $U$ so that
$F(z)=g_{z_0}(\Dt((z+i\be)/i))$, $z\in U$, where
$g_{z_0}(w)=w^{-\al}\,(1+\log w)^{-\frac1p}$ is a function of 1
complex variable with a determination of the log depending on
$z_0$. We remark that functions corresponding to two points
$z_0,z_1$ will only differ by constants which are irrelevant for
our estimates below, and for this reason we shall drop the
subindex in $g$.

We can now compute $\Box F(z)$ using the formula \Be \label{bern}
\Box[g(\Dt(z/i))]=\frac{(Bg)(\Dt(z/i))}{\Dt(z/i)},\quad z\in\Tu,
\Ee where $B=b(w\frac d{dw})$ is the 1 variable differential
operator of degree $r$ given by the Bernstein polynomial
$b(\la)=(-1)^r\la(\la+\frac d2)\cdots(\la+(r-1)\frac d2)$. One can
verify the equality (\ref{bern}) directly, using the Taylor series
of $g$ and $\Box[\Dt^n(z/i)]=b(n)\Dt^{n-1}(z/i)$ (see, e.g.,
\cite[p. 142]{FK}). Thus, an easy computation of the derivatives
of $g(w)$ leads to the expression
\[
|\Box F(z)|\leq \,
C\,|\Dt((z+i\be)/i)|^{-(\al+1)}\,(1+\log\Dt(y+\be))^{-\frac1p},\quad
z=x+iy\in\Tu,\] where we have also used $|\Dt(u+iv)|\geq \Dt(v)$,
$v\in\O$. Now, we can apply Lemmas \ref{ial} and \ref{log2} to
estimate the integral \[ \int_\O\left(\int_\SRn|\Box
F(x+iy)|^p\,dx\right)^\frac
qp\,\Dt^{\nu+q-\fnr}(y)\,dy\Pquad\Dquad
\]
\Beas \Dquad &\leq & C\,\int_\O\left(\int_\SRn
 \frac
{dx}{|\Dt(x+i(y+\be))|^{(\al+1)p}}\right)^\frac
qp\,\frac{\Dt^{\nu+q-\fnr}(y)\,dy}{(1+\log\Dt(y+\be))^\frac qp} \\
 &
\leq &
C'\,\int_\O\frac{\Dt^{\nu+q-\fnr}(y)}{\Dt(y+\be)^{((\al+1)p-\fnr)\frac
qp}}\,\frac{dy}{(1+\log\Dt(y+\be))^\frac qp}, \Eeas and observe
that the last quantity is finite for every $q\geq
p(1+\frac{\nu}{\fnr-1})$.
\smallskip

Let us now pass to the second type of counter-examples, and obtain
case (b) in the proposition. We may assume $q>2$. We know from
Proposition \ref{equivalent} that the inequality (\ref{cannot})
implies
 the existence of a constant
$C$ such that \Be\label{lpb}\|\sum f_j\|_p^q\leq C\sum
\Dt(\xi_j)^{-\nu}\| f_j\|_p^q,\Ee for any finite sequence
$\{f_j\}$ of Schwartz functions satisfying $\supp\hf_j\subset
B_{\frac 12}(\xi_j)$ and restricted to indices $j$ so that
$|\xi_j|<1$. Let us prove that this implies very easily the
necessity of $q<2q_\nu$ (which was already announced in
\cite{BBPR} for light-cones). Indeed, let us take
$f_j=\e_ja_je^{i(\xi_j|\cdot)}f$, with $\{\e_j\}$ a sequence of
Rademacher functions, and the support of $\hf$ a small
neighborhood of $0$. Taking the mean over all $\e_j$'s and using
Khintchine inequalities, we find that
$$\Bigl[\,\sum_j|a_j|^2\,\Bigr]^{\frac 12} \leq C \,\Bigl[\,\sum_j  \Dt(\xi_j)^{-\nu}|
a_j|^q\,\Bigr]^{\frac 1q},$$ with perhaps a different constant
$C$, independent of the sequence $\{a_j\}$. Now, choosing
$a_j=\Dt(\xi_j)^{\frac\nu{q-2}}$ and using $q>2$, this implies
that
$$\sum_{j;|\xi_j|<1}\Dt(\xi_j)^{\frac{2\nu}{q-2}}\leq C<\infty.$$
Using Proposition \ref{discr}, this
is equivalent to the fact that
$$\int_{|\xi|<1}\Dt(\xi)^{\frac{2\nu}{q-2}}\dxi<\infty,$$ which, in
turn, is equivalent to the condition $q<2q_\nu$. \ProofEnd

\bigskip


\section{Appendix: extension of the range on light-cones}
\setcounter{equation}{0} \setcounter{figure}{0}

In this additional section we just focus on  the previous problems
for the particular case of the light cone $\La_n$. As we shall
see, the sufficient conditions given in Proposition
\ref{equivalent} can be written in terms of Littlewood-Paley
inequalities in $\SR^n$, removing any dependence on complex
coordinates. One of such inequalities will lead us to a variant of
the so-called ``cone multiplier problem'', for which recent
results of Tao-Vargas and T. Wolff will provide us with some
positive answers to our question. We observe here that our
sufficient conditions are ``essentially necessary''. In fact,
Proposition \ref{equivalent} states that they are necessary for
the boundedness of the projector with other values of the
parameters. So negative results for them should also give new
regions of unboundedness for the Bergman projector. In this way,
our approach to this problem ends up in a series of challenging
questions in Harmonic Analysis, very closely related to the
present restriction and multiplier problems for the cone. One can
as well pose the corresponding questions for spheres, which even
in the 2 dimensional case seem to be unknown at present.

\medskip

Before particularizing for light-cones, we recall the situation
after the results in this paper. In Figure \ref{fig} we show the
regions of boundedness for the Bergman projector $P_\nu$ in $\Lpq$
spaces for a general symmetric cone. In the ``blank region'', our
main contribution up to now is Theorem \ref{necsuff}, which gives
the \emph{equivalence} between boundedness of $P_\nu$ in $\Lpq$
and the inequality \Be
\|\cLL\hf\|_\Lpq\lesssim \|f\|_{\Bpq},\quad f\in\Do. \label{R1}
\Ee   By self-adjointness of $P_\nu$ and the constraints explained
in the introduction, we only study this question in the ranges of
indices
\[
\nu>\tfrac n2-1, \quad 1\leq p<\infty, \quad 2\leq q<\tqn.
\]
We know the validity of (\ref{R1}) for all
$2<q<\min\{p,p'\}q_\nu$, and we know the sharpness of the right
index when $1\leq p\leq2$, so we shall also restrict our study to
$p>2$. In this case, and after the results in $\S4.4$,  the region
of uncertainty for the validity of (\ref{R1}) reduces to
\[
p'q_\nu\leq q<\min\{2q_\nu,\tqn\}.\] Our main contribution here is
to show that one can conclude positively for a small neighborhood around the
left endpoint of this interval.

\subsection{A Whitney covering for the light-cone}

From now on, we let $\O=\Lambda_n$ denote the light-cone in
$\SR^n$ with $n\geq3$. We first describe an explicit Whitney
decomposition for $\Lambda_n$.

\medskip

A natural candidate for a lattice in the light-cone  is
constructed as follows. For every $j\geq1$, take a maximal
$2^{-j}$-separated sequence $\{\om^{(j)}_k\}_{k=1}^{k_j}$ in the
sphere $S^{n-2}\subset\SR^{n-1}$, with respect to the Euclidean
distance (so that $k_j\sim 2^{j(n-2)}$). Then, define the
following grid of points in $\La_n$:\Be
\xi^\ell_{j,k}\;=\;\Bigl(2^\ell,2^\ell\,\sqrt{1-2^{-2j}}\,\om^{(j)}_{k}\Bigr),
\quad\ell\inZ,\; j\geq1,\; k=1,\ldots,k_j,\label{xiljk} \Ee and
the corresponding sets \Beas E^\ell_{j,k} & =&
\left\{(\tau,\xi')\in\Lambda_n\mid 2^{\ell-1}<\tau<2^{\ell+1},\;
{2^{-2j-2}}<1-\tfrac{|\xi'|^2}{\tau^2}<2^{-2j+2},\right.\Tquad\\
& &\left.
\Dquad\mand\left|\tfrac{\xi'}{|\xi'|}-\om^{(j)}_k\right|\leq
\dt2^{-j}\,\right\},
 \Eeas
 where the constant $\dt$ is chosen in such a way that the regions cover the cone. The geometric
picture in $\SR^3$ is as follows: the sets $E^\ell_{j,k}$ are
truncated conical shells  of height $\sim 2^\ell$, of thickness
$\sim 2^{\ell-2j}$, and further decomposed into $k_j\sim2^j$
sectors, all of equal arc-length $\sim 2^{\ell-j}$. This is the
usual decomposition of $\La_n$ for the study of cone multipliers
(see \cite{mock}). The next proposition proves that the regions
$E^\ell_{j,k}$ are very close from a Whitney covering of the cone.

\begin{proposition}
\label{Whit1} With the notation above, the grid $\{\xiljk\}$ is a
lattice in $\La_n$. Moreover, there exist $0<\eta_1<\eta_2$ such
that the corresponding family of invariant balls satisfy
\Benu
\item[(a)] $\{B_{\eta_1}(\xiljk)\}$ is disjoint in $\O$;
\item[(b)] $\{B_{\eta_2}(\xiljk)\}$ is a covering of $\O$;
\item[(c)] $B_{\eta_1}(\xiljk)\subset E^\ell_{j,k}\subset
B_{\eta_2}(\xiljk)$. \Eenu
\end{proposition}

\Proof

In view of the definition in $\S2.2$, it suffices to find two
fixed (invariant) balls $B\subset \tB$, centered at $\be$ and such
that \Be g^\ell_{j,k}(B)\,\subset\,E^\ell_{j,k}\,\subset
g^\ell_{j,k}(\tB),\label{inc1}\Ee for some fixed automorphisms of
the cone $g^\ell_{j,k}$ mapping $\be$ into $\xiljk$. Using
dilations as well as rotations with axis $\be$, it is sufficient
to prove this for the points $(1, \sqrt{1-2^{-2j}},0)=g_j\be$.
Then, an elementary exercise shows  that the corresponding set
$E_j$ is such that
$$E_j(c)\subset E_j\subset E_j(\tc),$$ where $E_j(c)$ is the
 set of all $\xi$ for which
$$ c^{-1}<\frac{(\xi|\be)}{(g_j\be|\be)}<c\;,\qquad c^{-1}<\frac
{\Dt(\xi)}{\Dt(g_j\be)}<c\;,\qquad c^{-1}<\frac
{\Dt_1(\xi)}{\Dt_1(g_j\be)}<c\;.$$ Here $\Dt_1(\xi)=
\xi_1-\xi_2$, and the constants $c,\tc$ are independent of $j$.
Let us finally prove that $g_j^{-1}(E_j(c))$ is contained in a
ball and contains a ball, centered at $\be$ and with radii
independent of $j$. From the invariance properties of the
quantities involved, this last set consists of $\xi$'s for which
$$c^{-1}<(g_j\xi|\be)<c\;,\qquad c^{-1}<
{\Dt(\xi)}<c\;,\qquad c^{-1}<{\Dt_1(\xi)}<c\;.$$ Using the
explicit value $(g_j\xi|\be)=\xi_1+\sqrt{1-2^{-2j}}\,\xi_2$, it is
an elementary exercise to find two such balls with radii
independent of $j$. \ProofEnd

\medskip

From this choice of points, and from the second condition in
Proposition \ref{equivalent}, we get the following result. To
simplify notation, we write $E_{j,k}=E^0_{j,k}$.

\begin{proposition}
{\bf : Weak sufficient condition.} \label{prR2}

Let $1\leq p,s<\infty$. Suppose that for some $\mu>0$ there exists
 $C_\mu$ such that  \Be \|\sum_{k=0}^{k_j}
f_k\|_p\leq\,C_\mu\, 2^{\frac {2\mu j}
s}\,\Bigl[\,\sum_{k=0}^{k_j}
 \|f_{k}\|^s_p\,\Bigr]^\frac1s,\quad
\forall\;j\geq1, \label{R23} \Ee for every sequence $\{f_{k}\}$
satisfying $\supp\hf_{ k}\subset E_{j,k}$. Then $P_\nu$ is bounded
in $\Lpq$ for all $q$ such that \Be \frac qs
<\min\left\{\frac{q_\nu}{q_\mu}, \frac{\nu-(\frac
n2-1)}\mu\right\}. \label{R66} \Ee
\end{proposition}

\Proof Using  Proposition \ref{equivalent}, we are reduced to
prove that the assumption implies the inequality $$
\|\sum_j\sum_{k=0}^{k_j} f_{j,k}\|_p\leq\,C_{\mu'}\,
\,\Bigl[\,\sum_j 2^{2j\frac {\mu'} s}\sum_{k=0}^{k_j}
 \|f_{j,k}\|^s_p\,\Bigr]^\frac1s,$$
with $\mu'$ perhaps larger, but arbitrarily close from $\mu$. Here
the Fourier transforms of the functions $f_{j, k}$ are contained
in the regions $E_{j,k}$. To prove the previous inequality, start
using Minkowski's inequality in $j$, apply (\ref{R66}) in each
block with fixed $j$, and conclude with H\"{o}lder's inequality.
\ProofEnd

\smallskip

\BR
  The natural conjecture in order to fill the whole ``blank
region'' in Figure \ref{fig} is that (\ref{R23}) holds, with $s=2$
and any $\mu>0$, within the range $2<p<\tfrac{2n}{n-2}$. This
range of $p$ coincides with the conjecture for the $\La(p)$-set
problem of $\SZ^n$-points in spheres (see \cite[5.5]{bourg2}). In
particular, it contains the conjectured range for the cone
multiplier in $\SR^n$ and for Bochner-Riesz in $\SR^{n-1}$:
$2<p<\tfrac{2(n-1)}{n-2}$. When $n=3$ the latter right end-point
is $p=4$, for which almost orthogonality techniques can be applied
to obtain some partial results (see next subsection). We observe
finally that, by Proposition \ref{equivalent}, our condition
(\ref{R23}) is necessary for the boundedness of $P_\mu$ in
$L^{p,s}_\mu$ when $\mu>\ind$.
 \ER

\subsection{Restriction techniques and new results for light-cones}

From now on we restrict to $p>2$. We shall study almost
orthogonality and restriction theorem techniques that can imply
our sufficient condition  \Be \|\sum_{k=0}^{k_j}
f_k\|_p\leq\,C_\e\, 2^{j\e}\,\Bigl[\,\sum_{k=0}^{k_j}
 \|f_k\|^2_p\,\Bigr]^\frac12, \label{R25} \Ee for $\e$ as close to $0$ as possible,
and the Fourier transforms of $f_k$ supported in $E_{jk}$. A
simple application of Minkowski's inequality shows that
(\ref{R25}) is implied by the \emph{square function estimate}\Be
\|\sum_{k=0}^{k_j} f_k\|_p\leq\,C_\e\,
2^{j\e}\,\Bigl\|\,\bigl(\sum_{k=0}^{k_j}
 |f_k|^2\,\bigr)^\frac12\,\Bigr\|_p\,. \label{R26} \Ee
 When $n=3$ and $p=4$ this has been
widely studied in relation with the cone multiplier problem. The
analogous question for the 2-dimensional disk (that is, for a
horizontal section of the cone) has a well-known positive answer
with $\e=0$, following from a simple geometric argument due to
C\'ordoba and Fefferman  \cite{cor}. This same argument is known
to be less sharp in the 3-dimensional cone,
where additional overlapping leads only to $\e=\frac14$ (see \cite{mock}).
Such estimate does not produce new results on the boundedness of
Bergman projectors, as one can check easily through
the numerology in our previous subsection.
One must beat the exponent $\e=\frac14$ to
obtain some improvement in our problem.

In this direction there are more recent works by
Bourgain \cite{bourg} and Tao-Vargas \cite{TV}, where this exponent has been
lowered to $\e=\frac14-\tau$, for a small $\tau>0$.
This improvement makes use of the so-called
\emph{bilinear restriction estimates} $R^*(2\times2\to q)$.
That is, finding the smallest value of $q\leq2$ for which: \Be
\bigl\|\,\widehat{g_1d\sigma_1}\;\widehat{g_2d\sigma_2}\,\bigr\|_{q}\lesssim
\|g_1\|_2\;\|g_2\|_2,\label{BR}\Ee when $g_1,g_2$ are smooth
functions with $1$-separated supports and $d\sigma$ is the surface
measure of the truncated cone
$\Lambda'_3=\{(|\xi'|,\xi')\in\SR^3\mid 1<|\xi'|<2\}$. More
precisely, we cite the following theorem from \cite{TV}, adapted
to our notation.

\begin{theorem}\label{TV} {\bf : see Theorem 5.1 in \cite{TV}.}
If $\kappa>0$ is such that $R^*(2\times2\to2-\kappa)$ holds for all
$g_1$, $g_2$ with unit separated support in $\Lambda'_3$, then for
all $\tau<\kappa/(16-4\kappa)$ the inequality (\ref{R26}) holds
with $\e=\frac14-\tau$.
\end{theorem}

Tao and Vargas are able to go down to $q>2-\frac8{121}$ in
(\ref{BR}), obtaining in (\ref{R26}) all $\e>\frac14-\frac1{238}$.
The recent sharp results for (\ref{BR}) given by T. Wolff in
\cite{W}, valid for all $q>2-\frac13$, improve the estimate in
(\ref{R26}) to all $\e>\frac14-\frac1{44}$. Whether one can go down
this index in (\ref{R26}) seems to be an open question, in spite
of the sharpness of Wolff's theorem in cone restriction. From this
discussion and a straighforward computation of the numerology we
conclude the following

\begin{corollary}
In $\SR^3$ and for $p=4$, the Bergman projector $P_\nu$ is bounded
in $L^{4,q}_\nu$ for all $2\leq q<
(\frac43+\frac1{24})\,q_\nu$, whenever $\nu>\frac12+\frac{5}{11}$.
\end{corollary}

By interpolation, this allows to give a positive answer in a part
of the blank region of Figure \ref{fig}. We do not give a precise
description of this new region since it will certainly be improved
in the future.

\medskip

 We conclude by observing that in
higher dimensions $n>3$, the sharp index for bilinear restriction
due to Wolff is $q>2-\frac{n-2}n$. We do not know if this has any
implication for the inequality (\ref{R26}), since the above cited
result of Tao-Vargas is only 3-dimensional.

\bigskip

\bibliographystyle{plain}

\end{document}